\documentclass[reqno,12pt,A4paper]{amsart}
\usepackage{graphicx}
\usepackage{epsfig}
\usepackage{amssymb,amsthm,amsmath,amsfonts,amstext}
\usepackage{latexsym}
\usepackage[all]{xy}

\newtheorem{thm}{Theorem}

\newtheorem{lem}{Lemma}
\newtheorem{prop}{Proposition}
\newtheorem{defin}{Definition}
\newcommand{\ignore}[1]{}

\input xy
\xyoption{all}

\def\d{\,{\rm{d}}}

\def\zm{\,{\zeta_{\mathcal{M}}  }  }
\def\proof{\,{\it Proof. }}

\title[Question mark function]
{Generating and zeta functions, structure, spectral and
analytic properties of the moments of Minkowski question mark function}
\author[Giedrius Alkauskas]{Giedrius Alkauskas}
\date{\today}

\begin{document}
\begin{abstract}
In this paper we are interested in moments of Minkowski question mark function $?(x)$. It appears that, to certain extent,  the results are analogous to the results obtained for objects associated with Maass wave forms: period functions, $L$-series, distributions, spectral properties. These objects can be naturally defined for $?(x)$ as well. Various previous investigations of $?(x)$ are mainly motivated from the perspective of metric number theory, Hausdorff dimension, singularity and generalizations. In this work it is shown that analytic and spectral properties of various integral transforms of $?(x)$ do reveal significant information about the question mark function. We prove asymptotic and structural results about the moments, calculate certain integrals involving $?(x)$, define an associated zeta function, generating functions, Fourier series, and establish intrinsic relations among these objects.
\end{abstract}
\maketitle
\begin{center}
\end{center}
\section{Introduction}
The aim of this paper to continue investigations on the moments of
Minkowski $F(x)$ function, begun in \cite{as} and \cite{ga2}.
The function $F(x)$ (``the question mark function") was introduced
by Minkowski in 1904 as an example of a continuous function
$F:[0,\infty)\rightarrow[0,1)$, which maps rationals to dyadic
rationals and quadratic irrationals to non-dyadic rationals. For
non-negative real $x$ it is defined by the expression
\begin{eqnarray}
F([a_{0},a_{1},a_{2},a_{3},...])=1-2^{-a_{0}}+2^{-(a_{0}+a_{1})}-2^{-(a_{0}+a_{1}+a_{2})}+...,
\label{min}
\end{eqnarray}
where $x=[a_{0},a_{1},a_{2},a_{3},...]$ stands for the
representation of $x$ by a (regular) continued fraction
\cite{khin}. More often this function is investigated in the interval $[0,1]$; in this case we use a standard notation $?(x)=2F(x)$ for $x\in[0,1]$. For
rational $x$, the series terminates at the last nonzero element
$a_{n}$ of the continued fraction. This function was investigated
by many authors. In particular, Denjoy \cite{denjoy} showed that
$?(x)$ is singular, and the derivative vanishes almost everywhere.
The proof of this fact is also given in \cite{as}. The problem arose in connection with the Calkin-Wilf tree (see below); the Farey tree is actually a subtree of it.
In fact, singularity of $?(x)$ follows from Khinchin's average value
theorem on continued fractions (\cite{khin}, chapter III). The
nature of singularity of $F(x)$ was clarified by Viader,
Parad\'{i}s and Bibiloni \cite{paradis1}. In particular, the
existence of the derivative $?'(x)$ in $\mathbb{R}$ for fixed $x$
forces it to vanish. Some other properties of $?(x)$ are
demonstrated in \cite{paradis2}. In the different direction, motivated
by the Hermite problem - to represent real cubic irrationals as periodic sequences of integers - Beaver and Garrity \cite{beaver} introduced a
two-dimensional analogue of $?(x)$. They showed that periodicity
of Farey iteration corresponds to a class of cubic irrationals,
and that the two dimensional analogue of $?(x)$ possesses similar
singularity properties. Nevertheless, the problem of Hermite still remains open. Salem in \cite{salem} proved (see also \cite{kinney}) that $?(x)$ satisfies the Lipschitz
condition of order $\frac{\log 2}{2\log \gamma}$, where $\gamma=\frac{1+\sqrt{5}}{2}$, and this is in fact the best possible exponent for Lipschitz condition. The Fourier-Stieltjes coefficients of $?(x)$, defined as $\int_{0}^{1}e^{2\pi i nx}\d ?(x)$, where investigated in \cite{salem}. The author, as an application of Wiener's theorem about Fourier series, gives average results on these coefficients without giving an answer to yet unsolved problem whether these coefficients vanish, as $n\rightarrow\infty$. It is worth noting that in Section 8 we will encounter analogous coefficients (see Proposition 3).\\
\indent Recently, Calkin and Wilf \cite{cw} defined a binary tree which
is generated by the iteration
$$
{a\over b}\quad\mapsto\quad {a\over a+b}\ ,\quad {a+b\over b},
$$
starting from the root ${1\over 1}$. Elementary considerations
show that this tree contains every positive rational number once and
only once, each being represented as a reduced fraction \cite{cw}.
The first four iterations lead to
\[
\xymatrix @R=.5pc @C=.5pc { & & & & & & & {1\over 1} & & & & & & & \\
& & & {1\over 2} \ar@{-}[urrrr] & & & & & & & & {2\over 1} \ar@{-}[ullll] & & & \\
& {1\over 3} \ar@{-}[urr] & & & & {3\over 2}\ar@{-}[ull] & & & & {2\over 3}\ar@{-}[urr] & & & & {3\over 1} \ar@{-}[ull] & \\
{1\over 4} \ar@{-}[ur] & & {4\over 3} \ar@{-}[ul] & & {3\over 5}
\ar@{-}[ur] & & {5\over 2} \ar@{-}[ul] & & {2\over 5} \ar@{-}[ur]
& & {5\over 3} \ar@{-}[ul] & & {3\over 4} \ar@{-}[ur] & & {4\over
1} \ar@{-}[ul] }
\]
Thus, the $n$th generation consists of $2^{n-1}$ positive
rationals. It is surprising that the
iteration discovered by M. Newman \cite{new}
$$
x_1=1,\quad x_{n+1}=1/(2[x_n]+1-x_n),
$$
produces exactly rationals of this tree, reading them line-by-line,
and thus gives an example of a simple recurrence which produces all
positive rationals. Recently, the authors of \cite{dilcher} produced a natural analogue of this tree, replacing integers $r$ with polynomials $r\in(\mathbb{Z}/2\mathbb{Z})[x]$. One of the results is that these polynomials also satisfy analogous recurrence (minding a proper definition of integral part of rational function, which comes for the Euclidean algorithm). It is important to note that the
$n$th generation of the Calkin-Wilf binary tree consists of exactly those
rational numbers, whose elements of the continued fraction sum up
to $n$. This fact can be easily inherited directly from the definition. First,
if rational number $\frac{a}{b}$ is represented as a continued fraction $[a_{0},a_{1},...,a_{r}]$, then the map $\frac{a}{b}\rightarrow\frac{a+b}{b}$ maps $\frac{a}{b}$ to $[a_{0}+1,a_{1}...,a_{r}]$. Second, the map $\frac{a}{b}\rightarrow\frac{a}{a+b}$ maps $\frac{a}{b}$ to $[0,a_{1}+1,...,a_{r}]$ in case $\frac{a}{b}<1$, and to $[1,a_{0},a_{1},...,a_{r}]$ in case $\frac{a}{b}>1$. This is an important fact which makes the investigations of rational numbers according to their position in the Calkin-Wilf tree highly motivated from the perspective of metric number theory and dynamics of continued fractions. The sequence of numerators
\begin{eqnarray*}
0,1,1,2,1,3,2,3,1,4,3,5,2,5,3,4,1,...
\end{eqnarray*}
is called the Stern diatomic sequence and was introduced in \cite{stern}. It satisfies the recurrence relations
\begin{eqnarray*}
s(0)=0,\quad s(1)=1,\quad s(2n)=s(n),\quad
s(2n+1)=s(n)+s(n+1).
\end{eqnarray*}
This sequence and the pairs $(s(n),s(n+1))$ have also
investigated by Reznick \cite{reznick}. It is not surprising (minding the relation to Farey tree) that the ``distribution" of numerators, which are defined via the moments $Q^{(\tau)}_{N}
=\sum_{n=2^{N}+1}^{2^{N+1}}s^{2\tau}(n)$, for $\tau>0$, has an interesting application in thermodynamics and spin physics \cite{cvit}, \cite{cont}.\\
\indent In \cite{as} it was shown that each generation of the
Calkin-Wilf tree possesses a distribution function $F_{n}(x)$, and that
$F_{n}(x)$ converges uniformly to $F(x)$. This is, of course, a well known fact about the Farey tree. The function $F(x)$ as a
distribution function is uniquely determined by the functional
equation \cite{as}
\begin{eqnarray}
2F(x)=\left\{\begin{array}{c@{\qquad}l} F(x-1)+1 & \mbox{if}\quad
x\geq 1,
\\ F({x\over 1-x}) & \mbox{if}\quad 0\leq x<1. \end{array}\right.
\label{distr}
\end{eqnarray}
This implies $F(x)+F(1/x)=1$. The mean value of $F(x)$ has been
investigated by several authors, and was proved to be $3/2$
(\cite{as}, \cite{reznick}).\\
\indent On the other hand, almost all the results mentioned reveal the properties of the Minkowski question mark function as function itself. Nevertheless, the final goal and motivation of papers \cite{as}, \cite{ga2} and this work is to show that in fact there exist several unique and very interesting analytic objects associated with $F(x)$ which encode a great deal of essential information about it. These objects will be introduced in Section 2.
Lastly, and most importantly, let us point out that, surprisingly, there are
striking similarities and analogies between the results proved here as well as in \cite{ga2}, with Lewis-Zagier's \cite{zagl} results on
period functions for Maass wave forms. That work is an expanded and clarified exposition of an earlier paper by Lewis \cite{lewis}. The concise exposition of these objects, their properties and relations to Selberg zeta function can be found in \cite{zagier}. The
reader who is not indifferent to the beauty of Minkowski's question mark
function is strongly urged to compare results in this work with those in
\cite{zagl}. Thus, instead of making quite numerous
references to \cite{zagl} at various stages of the work (mainly in Sections 2, 3, 8 and 9), it is more useful to
give a table of most important functions encountered there, juxtaposed with
analogous object in this work. Here is the summary (the notations on the right
will be explained in Sections 2 and 9). \it\\

\begin{tabular}{|r | l || r| l|}
\hline
Maass wave form    & $u(z)$                    & $\Psi(x)$ & Periodic function on the real line\\
Period function    &  $\psi(z)$                & $G(z)$    & Dyadic period function\\
Distribution       & $U(x)\d x$                & $\d F(x)$ & Minkowski's "question mark"\\
$L-$functions      & $L_{0}(\rho),L_{1}(\rho)$ & $\zm(s)$  & Dyadic zeta function\\
Entire function    & $g(w)$                    & $m(t)$    & Generating function of moments\\
Entire function    & $\phi(w)$                 & $M(t)$    & Generating function of moments\\
Spectral parameter & $s$               & $\frac{1}{2}$; $1$ & Analogue of spectral parameter\\
\hline
\end{tabular}\\

\rm As a matter of fact, the first entry is the only one where
the analogy is not precise. Indeed, the distribution $U(x)$ is the limit value
of the Maass wave form $u(x+iy)$ on the real line (as $y\rightarrow +0$),
in the sense that $u(x+iy)\sim y^{1-s}U(x)+y^{s} U(x)$, whereas $\Psi(x)$ is the same $F(x)$ made periodic. As far as the last entry
of the table in concerned, the ``analogue" of spectral parameter,
sometimes this role is played by $1$, sometimes by $\frac{1}{2}$. This occurs,
obviously, because the relation between Maass form and $F(x)$ is only the
analogy which is not strictly defined.\\

This work is organized as follows. In Section 2 we give a summary of the previous results, obtained in \cite{as} and \cite{ga2}. In Section 3 we give a short proof of three term functional equation (\ref{sim}), and prove the existence of certain distributions, which can be thought as close relatives of $F(x)$. In Section 4 it is demonstrated that there are linear relations among moments $M_{L}$, and they are presented in an explicit manner. Moreover, we formulate a conjecture, based on the analogy with periods, that these are the only possible relations. In Section 5, the estimate for the moments $m_{L}$ is proved. As a consequence, $\lim_{L\rightarrow\infty}\frac{\log m_{L}}{\sqrt{L}}=-2\sqrt{\log{2}}$. In Section 6 we prove the exactness of a certain sequence of functional vector spaces and linear maps related to $F(x)$ in an essential way. Section 7 is devoted to calculation of number of integrals, giving a rare example of Stieltjes integral, involving the question mark function, that ``can" be calculated. In Section 8 we compute the Fourier expansion of $F(x)$. It is shown that this establishes yet another relation among $m(t)$, $G(z)$ and $F(x)$ via Taylor coefficients and special values. In the penultimate Section 9, the associated Dirichlet series $\zm(s)$ is introduced. In the last section, some concluding remarks are presented, regarding future research; relations between $F(x)$ and the Calkin-Wilf tree (and the Farey tree as well) to the known objects are established. Note also that we use the word ``distribution" to describe a monotone function on $[0,\infty)$ with variation $1$, and also for a continuous linear functional on some space of analytic functions. In each case the meaning should be clear from the context. \\
\begin{figure}[h]
\begin{center}
\includegraphics[width=210pt,height=430pt,angle=-90]{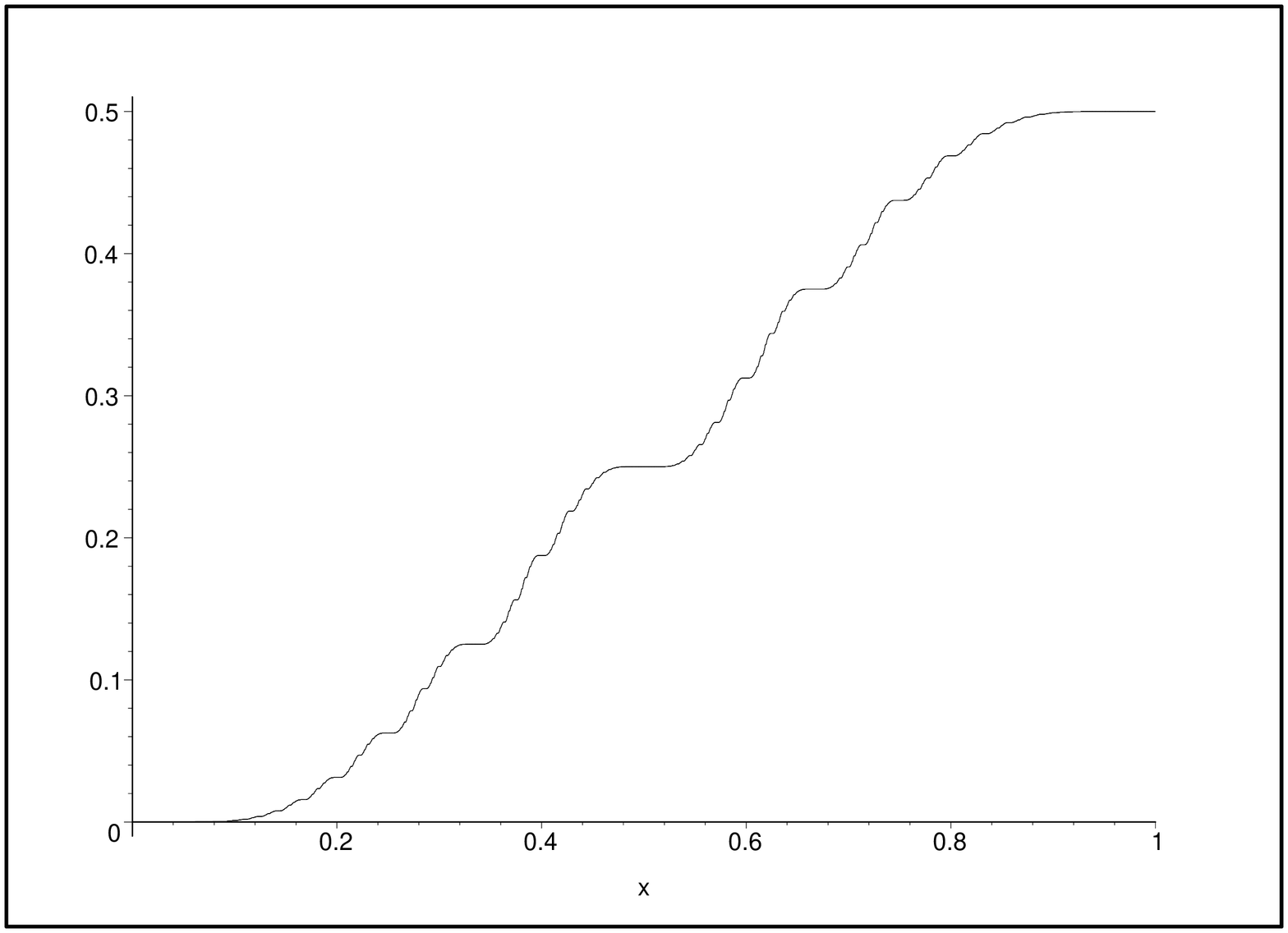}
\caption{Minkowski's question mark $?(x)$, $x\in[0,1]$}
\end{center}
\end{figure}
\section{Summary of previous results}
This section provides a summary of previous results. Let
\begin{eqnarray}
M_{L}=\int\limits_{0}^{\infty}x^{L}\d F(x),\quad
m_{L}=\int\limits_{0}^{\infty}\Big{(}\frac{x}{x+1}\Big{)}^{L}\d
F(x)=2\int\limits_{0}^{1}x^{L}\d F(x).\label{moments}
\end{eqnarray}
Both sequences are of definite number-theoretical significance because
\begin{eqnarray}
M_{L}=\lim_{n\rightarrow\infty}2^{1-n}\sum\limits_{a_{0}+a_{1}+...+a_{s}=n}[a_{0},a_{1},..,a_{s}]^{L},\quad
m_{L}=\lim_{n\rightarrow\infty}2^{2-n}\sum\limits_{a_{1}+...+a_{s}=n}[0,a_{1},..,a_{s}]^{L},
\label{mom}
\end{eqnarray}
(the summation takes place over rational numbers presented as continued fractions; thus, $a_{i}\geq 1$ and $a_{s}\geq 2$). In fact, clarification of their nature was the initial main motivation for our work.
We define the exponential generating functions
$M(t)=\sum_{L=0}^{\infty}\frac{M_{L}}{L!}t^{L}$,
$m(t)=\sum_{L=0}^{\infty}\frac{m_{L}}{L!}t^{L}$. Thus,
\begin{eqnarray*}
M(t)=\int\limits_{0}^{\infty}e^{xt}\d F(x),\quad m(t)=\int\limits_{0}^{\infty}\exp\Big{(}\frac{xt}{x+1}\Big{)}\d F(x)=2\int\limits_{0}^{1}e^{xt}\d F(x).
\end{eqnarray*}
One easily verifies that $m(t)$ is an
entire function and that $M(t)$ has radius of convergence $\log 2$. There
are natural relations among values $M_{L}$ and $m_{L}$,
independent of a specific distribution, like $F(x)$. They encode the
relations among functions $x^{L}$, $L\in\mathbb{N}_{0}$, and
functions $\big{(}\frac{x}{x+1}\big{)}^{L}$, $L\in\mathbb{N}_{0}$, given by
$x^{L}=\sum\limits_{s\geq L}\binom{s-1}{L-1}\big{(}\frac{x}{x+1}\big{)}^{s}$.
Therefore,
\begin{eqnarray}
M_{L}=\sum\limits_{s\geq L}\binom{s-1}{L-1}m_{s}.
\label{est}
\end{eqnarray}
On the other hand, the intrinsic information about $F(x)$ is
encoded in the relation
\begin{eqnarray}
m_{L}=M_{L}-\sum\limits_{s=0}^{L-1}M_{s}\binom{L}{s},\quad L\geq
0. \label{intr}
\end{eqnarray}
Further, we have
\begin{eqnarray}
M(t)=\frac{1}{2-e^{t}}m(t),\quad m(t)=e^{t}m(-t). \label{symm}
\end{eqnarray}
The first relation is equivalent to the system (\ref{intr}), and it encodes all the
information about $F(x)$ (provided we take into account the natural relations
just mentioned). The second one represents only the symmetry property, given by
$F(x)+F(1/x)=1$. One of the main results about $m(t)$ is that it is uniquely
determined by the regularity condition $m(-t)\ll e^{-\sqrt{t\log 2 }}$, as
$t\rightarrow\infty$, the boundary condition $m(0)=1$, and the integral
equation
\begin{eqnarray}
m(-s)=(2e^{s}-1)\int\limits_{0}^{\infty}m'(-t)J_{0}(2\sqrt{st})\d t,\quad
s\in\mathbb{R}_{+}. \label{integ}
\end{eqnarray}
(Here $J_{0}(*)$ stands for the Bessel function
$J_{0}(z)=\frac{1}{\pi}\int_{0}^{\pi}\cos(z\sin x)\d x$). This equation can be
rewritten as a second type Fredholm integral equation (\cite{kolm}, chapter 9).
In fact, if we denote
\begin{eqnarray*}
\psi(s)=\sqrt{2e^{s}-1},\quad
\frac{J_{1}(2\sqrt{st})}{\psi(s)\psi(t)}=K(s,t), \quad
\frac{m(-s)-1}{\sqrt{s}\psi(s)}={\sf m}(s),
\end{eqnarray*}
then one has
\begin{eqnarray}
{\sf m}(s)=\ell(s)-\int\limits_{0}^{\infty}{\sf m}(t)K(s,t)\d t,\text{ where
}\ell(s)=-\frac{1}{\psi(s)}\int\limits_{0}^{\infty}
\frac{J_{1}(2\sqrt{st})}{\sqrt{t}(2e^{t}-1)}\d t.\label{fred}
\end{eqnarray}
On the other hand, all the results about exponential generating
function can be restated in terms of generating function of
moments. Let $G(z)=\sum\limits_{L=0}^{\infty}m_{L+1}z^{L}$ for $|z|\leq 1$. Moreover, the functional equation for $G_{\lambda}(z)$ implies that there exist all left derivatives of $G_{\lambda}(z)$ at $z=1$. Then the integral
\begin{eqnarray}
G(z)=\int\limits_{0}^{\infty}\frac{\frac{x}{x+1}}{1-\frac{x}{x+1}z}\d F(x)=2\int\limits_{0}^{1}\frac{x}{1-xz}\d F(x).
\label{gen}
\end{eqnarray}
extends $G(z)$ to the cut plane $\mathbb{C}\setminus(1,\infty)$. The generating
function of moments $M_{L}$ does not exists due to the factorial growth of
$M_{L}$, but the generating function can still be defined in the cut plane
$\mathbb{C}'=\mathbb{C}\setminus(0,\infty)$ by
$\int_{0}^{\infty}\frac{x}{1-xz}\d F(x)$. In fact, this integral equals to
$G(z+1)$, which is the consequence of an algebraic identity
\begin{eqnarray*}
\frac{x}{1-xz}=\frac{\frac{x}{x+1}}{1-\frac{x}{x+1}(z+1)}.
\end{eqnarray*}
The following result was proved in \cite{ga2}.
\begin{thm}
The function $G(z)$, defined initially as a power series, has an analytic
continuation to the cut plane $\mathbb{C}\setminus(1,\infty)$ via (\ref{gen}). It
satisfies the functional equation
\begin{eqnarray}
-\frac{1}{1-z}-\frac{1}{(1-z)^{2}}G\Big{(}\frac{1}{1-z}\Big{)}+2G(z+1)=G(z),
\label{funct}
\end{eqnarray}
and also the symmetry property
\begin{eqnarray*}
G(z+1)=-\frac{1}{z^{2}}G\Big{(}\frac{1}{z}+1\Big{)}-\frac{1}{z}.
\end{eqnarray*}
Moreover, $G(z)\rightarrow 0$, if $z\rightarrow\infty$ and the distance from $z$ to a half line $[0,\infty)$ tends to infinity.\\
Conversely, the function having these properties is unique.
\end{thm}
Note that two functional
equations for $G(z)$ can be merged into single one. It is easy
to check that the equation
\begin{eqnarray}
\frac{1}{z}+\frac{1}{z^{2}}G\Big{(}\frac{1}{z}\Big{)}+2G(z+1)=G(z)\label{sim}
\end{eqnarray}
is equivalent to both of them together. In fact, the change $z\mapsto 1/z$
in the last equation gives the symmetry property, and application
of it to the term $G(1/z)$ in (\ref{sim}) gives the functional
equation in the Theorem 1. Nevertheless, it is sometimes convenient to separate (\ref{sim}) into two equations. The reason for this is that in (\ref{funct}) all arguments belong simultaneously to $\mathbb{H}$, $\mathbb{R}$, or $\mathbb{H}^{-}$, whereas in (\ref{sim}) they are mixed. This will become crucial later (see the last the section).\\
The transition $m(t)\rightarrow G(z)$ is given by Laplace transform:
\begin{eqnarray*}
1+zG(z)=\int\limits_{0}^{\infty}m(zt)e^{-t}\d t.
\end{eqnarray*}
The same transform applied to the eigen-functions of the Fredholm integral
equation (\ref{fred}), yields the following result \cite{ga2}.
\begin{thm}
For every eigen-value $\lambda$ of the integral operator
associated with the kernel $K(s,t)$, there exists at least one
holomorphic function $G_{\lambda}$ (defined for
$z\in\mathbb{C}\setminus\mathbb{R}_{>1}$), such that the
following holds:
\begin{eqnarray}
2G_{\lambda}(z+1)=G_{\lambda}(z)+\frac{1}{\lambda
z^{2}}G_{\lambda}\Big{(}\frac{1}{z}\Big{)}.\label{eigen}
\end{eqnarray}
Moreover, $G_{\lambda}(z)$ for $\Re z<0$ satisfies all regularity
conditions, imposed by it being an image under Laplace transform
(\cite{lavr}, p. 469).\\
Conversely: for every $\lambda$, such that there exists a
function, which satisfies $(\ref{eigen})$ and these conditions,
$\lambda$ is the eigen-value of this operator. The set of all
possible $\lambda$'s is countable, and $\lambda_{n}\rightarrow 0$,
as $n\rightarrow\infty$.
\end{thm}
\begin{figure}
\centering
\begin{tabular}{c c}

\epsfig{file=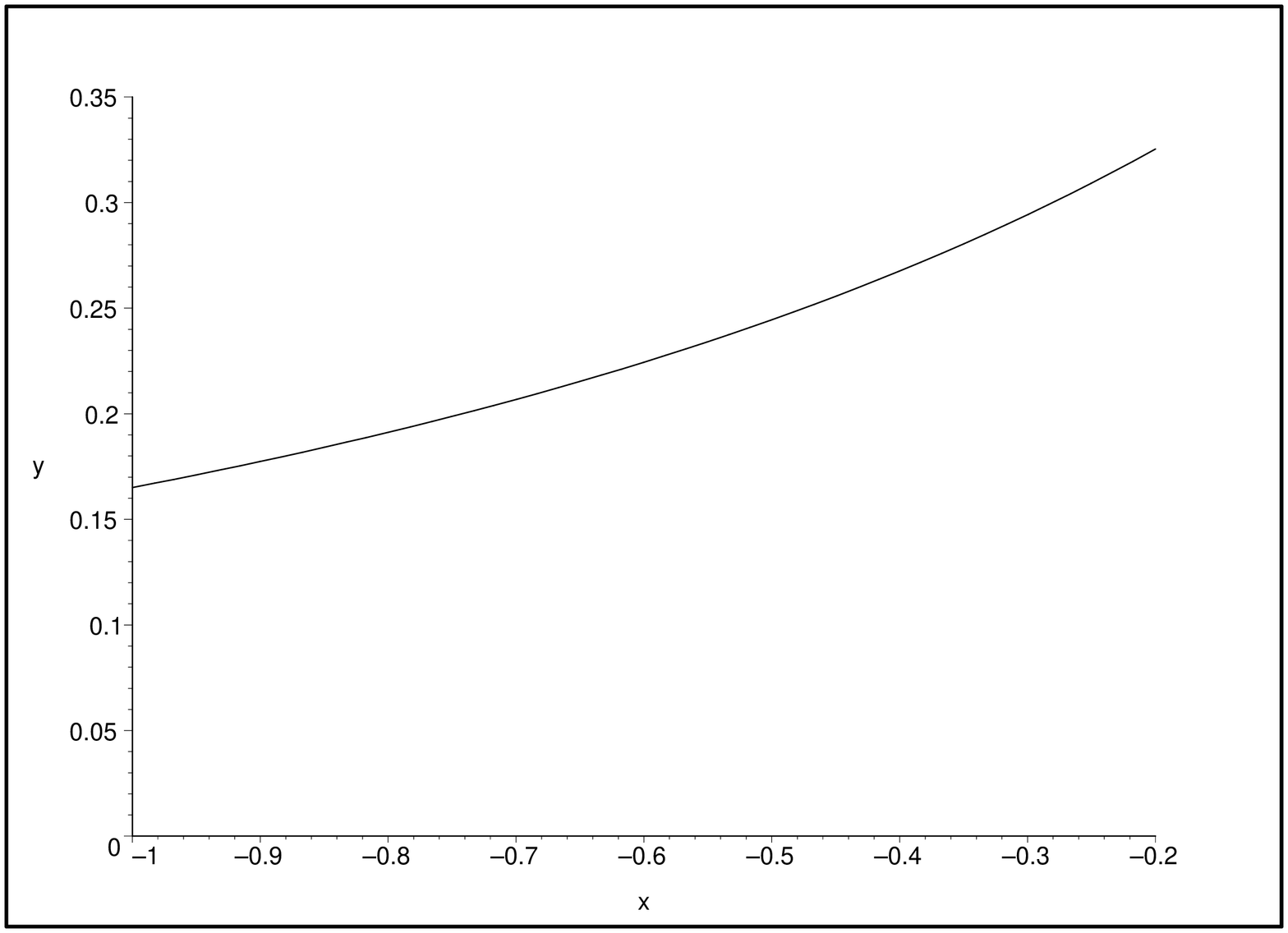,width=0.37\textwidth,angle=-90}
&\epsfig{file=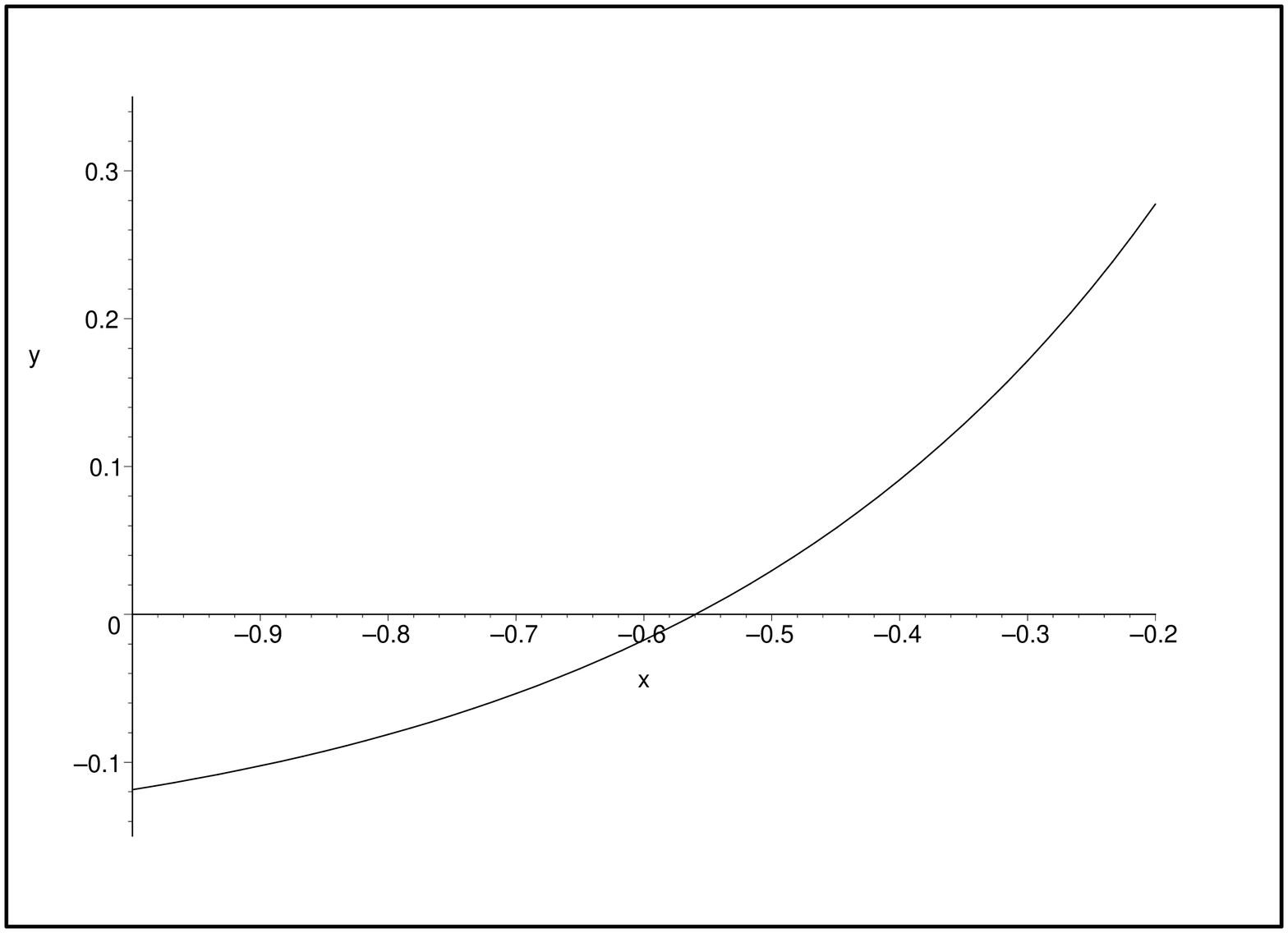,width=0.37\textwidth,angle=-90} \\
$\lambda_{1}=0.25553210$ & $\lambda_{2}=-0.08892666$\\

\epsfig{file=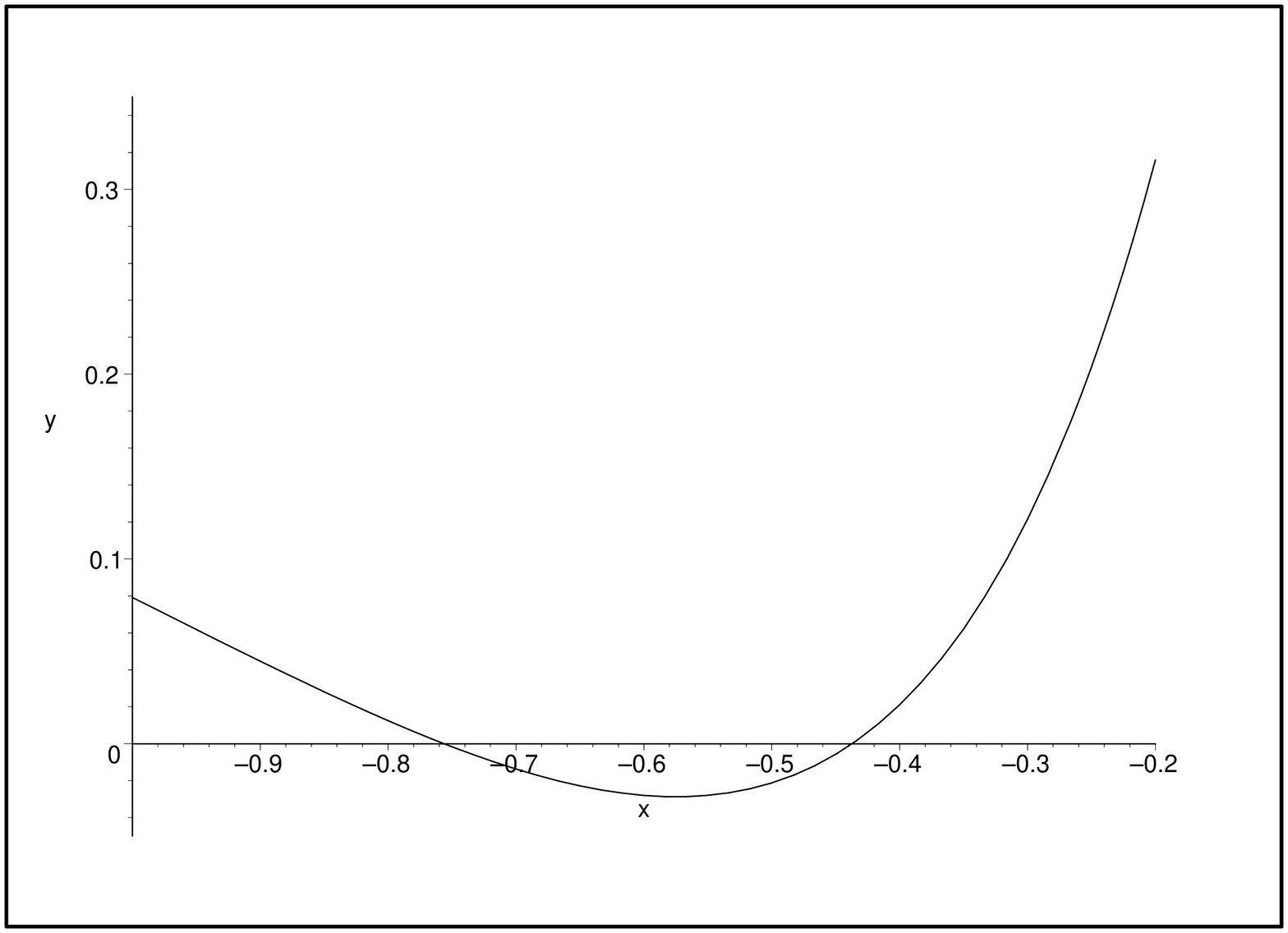,width=0.37\textwidth,angle=-90} &
\epsfig{file=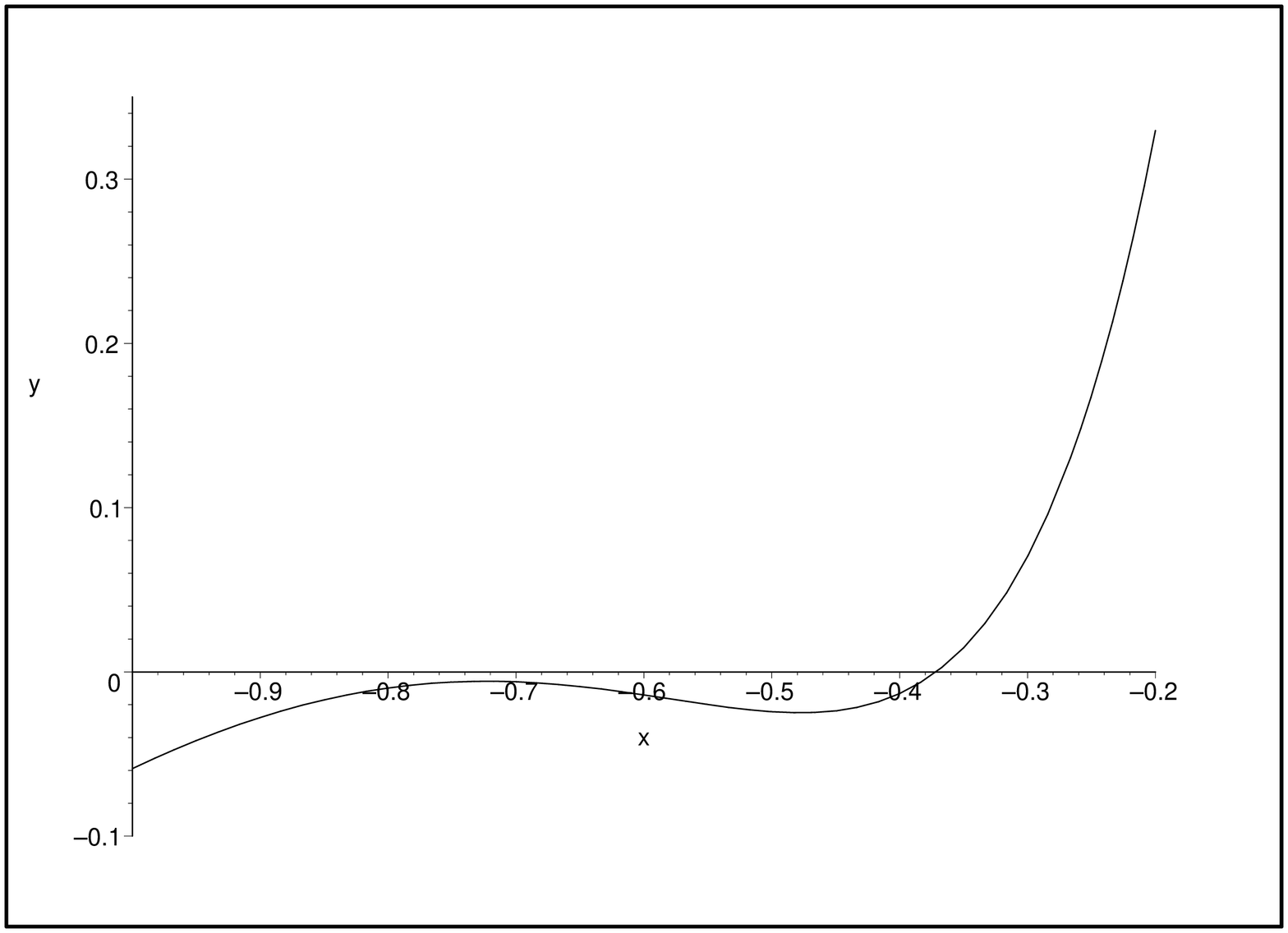,width=0.37\textwidth,angle=-90}\\
$\lambda_{3}=0.03261586$ & $\lambda_{4}=-0.01217621$\\

\epsfig{file=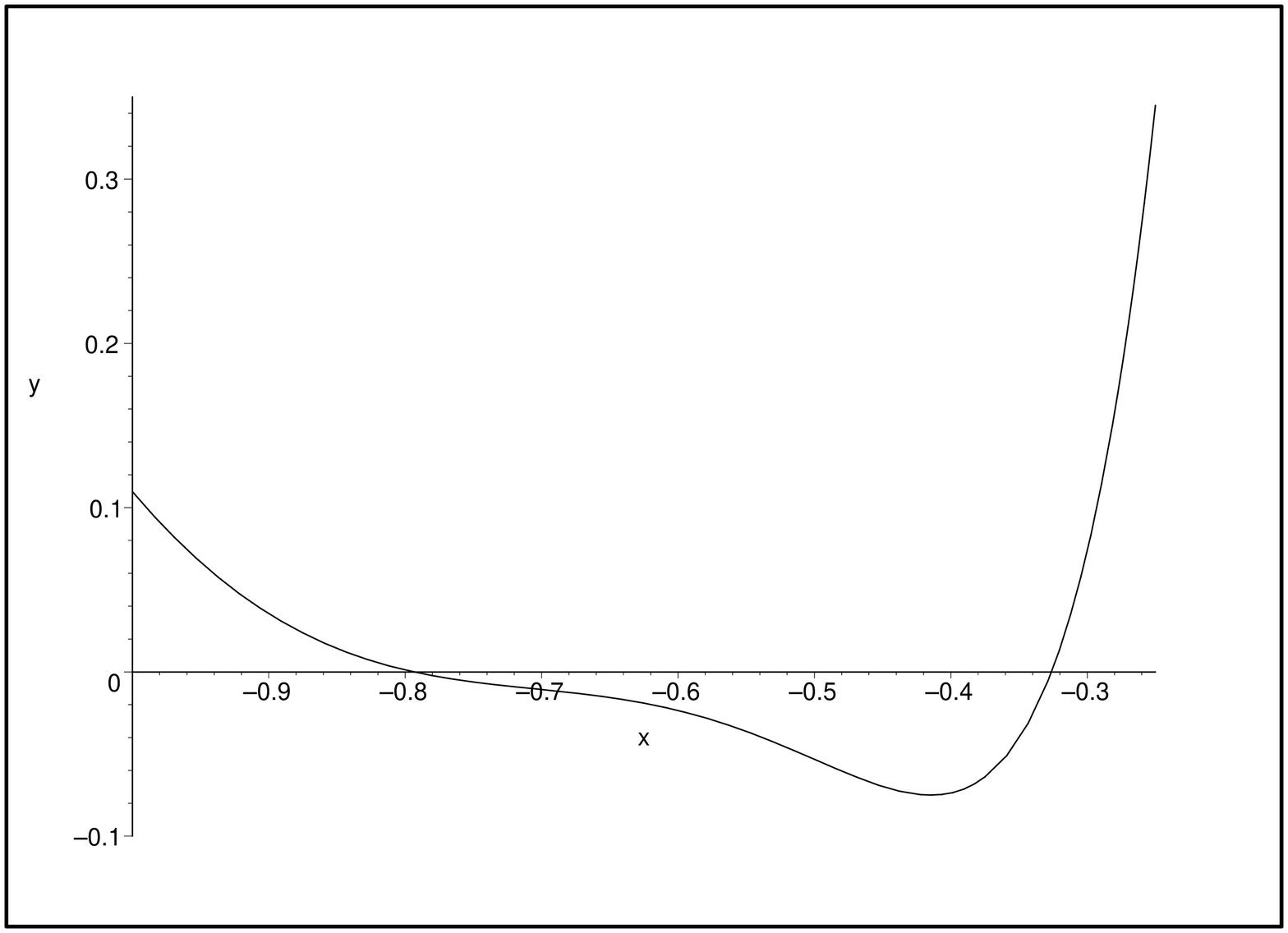,width=0.37\textwidth,angle=-90} &
\epsfig{file=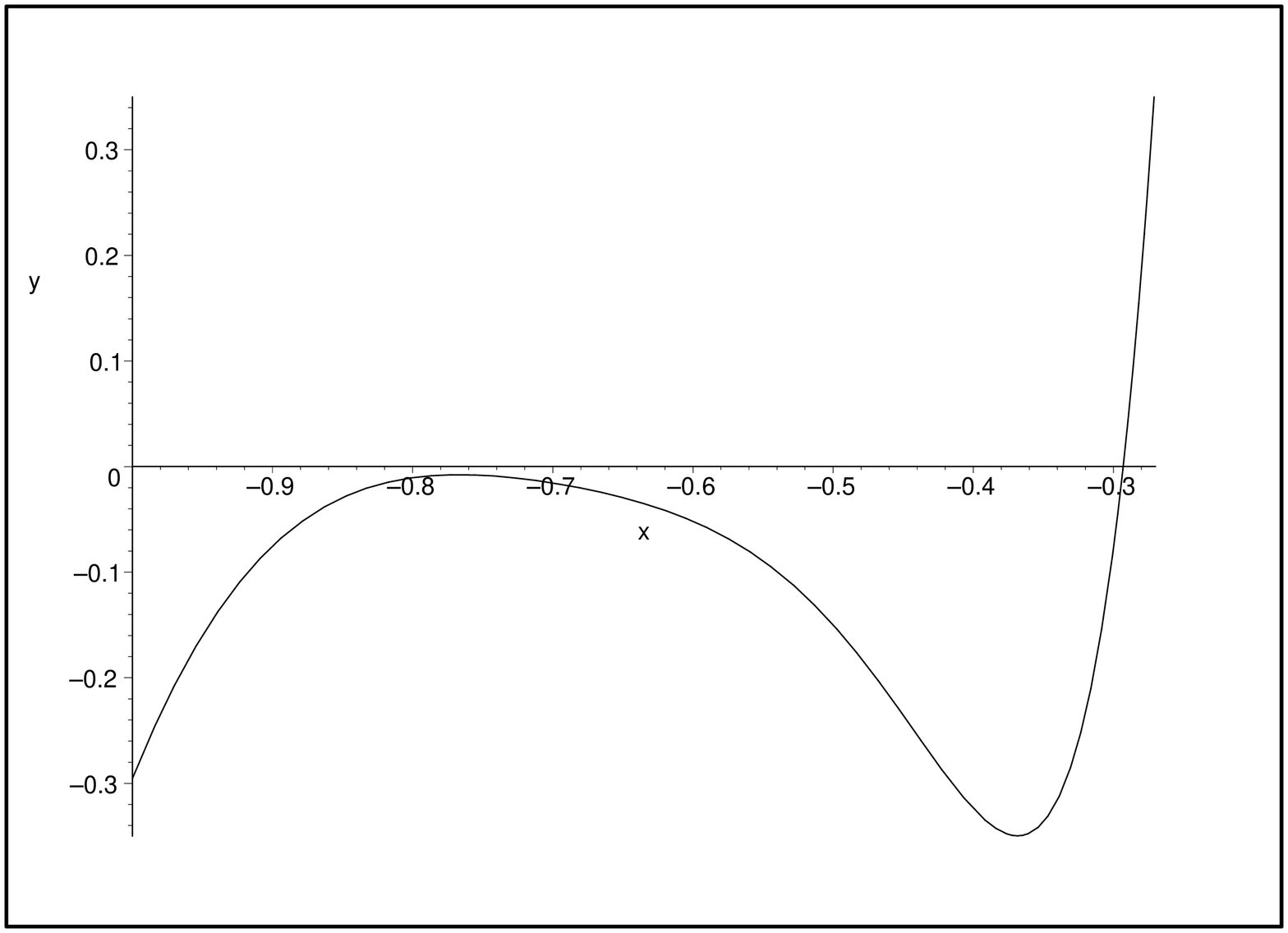,width=0.37\textwidth,angle=-90}\\
$\lambda_{5}=0.00458154$ & $\lambda_{6}=-0.00173113$\\
\end{tabular}
\caption{Eigen-functions $G_{\lambda}(z)$ for $z\in[-1,-0.2]$}.
\end{figure}
Figure 2 shows the functions $G_{\lambda}(z)$ (for the first six eigen-values) for real $z$ in the interval $[-1,-0.2]$. The choice of this interval is motivated by Theorem 5. Note also that functional equation implies $G_{\lambda}(0)=\big{(}\frac{1}{2}+\frac{1}{2\lambda}\big{)}G_{\lambda}(-1)$. Thus, one has $\frac{G_{\lambda}(0)}{G_{\lambda}(-1)}\rightarrow\infty$, as $\lambda\rightarrow0$. This can also be seen empirically from Figure 2.\\
Summarizing, there are three objects associated with the Minkowski question mark function.\it
\begin{itemize}
\item Distribution $F(x)$ = Functional equations (\ref{distr}) + Continuity.\\
\item Period function $G(z)$ = Three term functional equation (\ref{sim}) +
Mild
growth condition (as in Theorem 1).\\
\item Exponential generating function $m(t)$ = Integral equation (\ref{integ})
+
boundary value and diminishing condition on the negative real line.\\
\end{itemize}\rm
Each of them is characterized by the functional equation, and subject to some
regularity conditions, it is unique, and thus arises exactly from $F(x)$. The objects are described via the ``equality" {\it Function =
Equation + Condition}. This means that the object on the left possesses both
features; conversely - any object with these properties is necessarily the
function on the left.\\
As expected, here we encounter the phenomena of ``bootstrapping": in all cases,
regularity conditions can be significantly relaxed, and they are sufficient
for the uniqueness, which automatically implies stronger regularity conditions.
Here we show the rough picture of this phenomena. In each case, we suppose that
the object satisfies the corresponding functional equation. For the details,
see \cite{ga2}.
\begin{eqnarray*}
&(i)& F(x) \text{ is continuous at one point }\Rightarrow F(x) \text{ is continuous}.\\
&(ii)&\text{ There exists } \varepsilon<1 \text{ such that for every } z \text{
with }\Re z<0, \text{ we have }\\&\quad& G(z-x)=O(2^{\varepsilon x}) \text{ as
} x\rightarrow\infty\Rightarrow G(z)=
O(|z|^{-1}) \text{ as } dist(z,\mathbb{R}_{+})\rightarrow\infty.\\
&(iii)& m'(-t)=O(t^{-1}) \text{ as }t\rightarrow\infty\Rightarrow |m(-t)|\ll
e^{-\sqrt{t\log 2 }}\text{ as }t\rightarrow\infty.
\end{eqnarray*}
Corresponding converse results were proved in \cite{ga2}. For $F(x)$, this was
in fact the starting point of these investigations, since the distribution of
rationals in the Calkin-Wilf tree is a certain continuous function satisfying
(\ref{distr}); thus, it is exactly $F(x)$. The converse result for $m(t)$
follows from Fredholm alternative, since all eigen-values of the operator
(\ref{fred}) are strictly less than $1$ in absolute value. Finally, the
converse theorem for $G(z)$ follows from a technical detail in the proof, which is
the numerical estimate $0<\frac{\pi^{2}}{12}-\frac{\log^{2}2}{2}<1$; as a matter
of fact, it appears that this is essentially the same argument as in the case of $m(t)$, since this constant gives the upper bound for the moduli of
eigen-values.\\
One of the aims of this paper is to clarify the connections among these three objects,
and to add the final fourth satellite, associated with $F(x)$. Henceforth, we
have the complete list:\it\\
\begin{itemize}
\item Zeta function $\zm(s)$ (see definition (\ref{dir}) below) = Functional
equation with symmetry $s\rightarrow -s$ (\ref{tet}) + Regularity behavior in
vertical strips.
\end{itemize}
\rm
In this case, we do not present a proof of a converse result. Indeed, the converse result for $G(z)$ is strongly motivated by its relation to Eisenstein series $G_{1}(z)$ (see \cite{ga2} and the last Section). In the case of $\zm(s)$, this question is of small importance, and we rather concentrate on the direct result and its consequences.
\section{Three term functional equation, distributions $F_{\lambda}(x)$}
In this section, we give another proof of (\ref{sim}), different from the one presented in \cite{ga2}, since it is considerably shorter. For our purposes, it is convenient to
work in a slightly greater generality. Suppose that $\lambda\in\mathbb{R}$ has the property that there exists a function $F_{\lambda}(x)$, $x\in[0,\infty)$, such that
\begin{eqnarray}
 \d F_{\lambda}(x+1)=\frac{1}{2}\d F_{\lambda}(x), \d F_{\lambda}\Big{(}\frac{1}{x}\Big{)}=\frac{1}{\lambda}\d F_{\lambda}(x).\label{psi}
\end{eqnarray}
 We omitted the word ``continuous" in the description of the function intentionally. For a moment, consider $F_{\lambda}(x)=F(x)$ with $\lambda=-1$. Then $F_{-1}(x)$ is certainly continuous. The reason for introducing $\lambda$ will be apparent later.
Thus, let
\begin{eqnarray*}
G_{\lambda}(z)=\int\limits_{0}^{\infty}\frac{1}{x+1-z}\d F_{\lambda}(x).
\end{eqnarray*}
Since $F(x)+F(1/x)=1$, we see that for $\lambda=-1$ this agrees with the definition
(\ref{gen}). This converges to an analytic function in the cut plane
$\mathbb{C}\setminus(1,\infty)$. We have
\begin{eqnarray*}
2G_{\lambda}(z+1)=2\int\limits_{0}^{1}\frac{1}{x-z}\d F_{\lambda}(x)+
2\int\limits_{1}^{\infty}\frac{1}{x-z}\d F(x)=\\
2\int\limits_{0}^{\infty}\frac{1}{\frac{x}{x+1}-z}\d F_{\lambda}\Big{(}\frac{x}{x+1}\Big{)}+
2\int\limits_{0}^{\infty}\frac{1}{x+1-z}\d F_{\lambda}(x+1)=\\
\frac{2}{z}\int\limits_{0}^{\infty}
\Big{(}\frac{x+1}{x+1-\frac{1}{z}}-
1+1\Big{)}\d F_{\lambda}\Big{(}\frac{1}{x+1}\Big{)}+G_{\lambda}(z)=\\
\frac{\alpha}{\lambda z}+\frac{1}{\lambda z^{2}}
G_{\lambda}\Big{(}\frac{1}{z}\Big{)}+G_{\lambda}(z),
\text{ where } \alpha=\int\limits_{0}^{\infty}\d F_{\lambda}(x).
\end{eqnarray*}
For $\lambda=-1$ and $F_{-1}(x)=F(x)$,
this gives Theorem 1. Further, suppose $\lambda\neq -1$. Then
\begin{eqnarray*}
\alpha=\int\limits_{0}^{\infty}\d F_{\lambda}(x)=\int\limits_{1}^{\infty}\d F_{\lambda}(x)
+\int\limits_{0}^{1}\d F_{\lambda}(x)=\frac{\alpha}{2}-\frac{\alpha}{2\lambda}\Rightarrow \alpha=0.
\end{eqnarray*}
Therefore, the last functional equation reads as
\begin{eqnarray*}
2G_{\lambda}(z+1)=\frac{1}{\lambda z^{2}}G_{\lambda}\Big{(}\frac{1}{z}\Big{)}+G_{\lambda}(z).
\end{eqnarray*}
As a matter of fact, there cannot be any reasonable functions $F_{\lambda}(x)$, satisfying (\ref{psi}). Nevertheless, the last functional equation is identical to (\ref{eigen}). Thus, Theorem 2 gives a description of all such possible $\lambda$. This suggests that we can still find certain distributions $F_{\lambda}(x)$. Further, as it was mentioned, $-1$ is not an-eigen value of (\ref{fred}). Due to the minus sign in front of the operator, this is exactly the exceptional eigen-value, which is essential in the Fredholm alternative. The above proof (rigid at least in case $\lambda=-1$), surprisingly, proves that the next tautological sentence has a certain point: ``$-1$ is not an eigen-value because it is $-1$". Indeed, we obtain non-homogeneous part of the three term functional equation only because $\lambda=-1$, since otherwise $\alpha=0$ and the equation is homogenic.\\
Distributions $F_{\lambda}(x)$ can indeed be strictly defined, at least in the space ${\sf C}^{\omega}$ of functions, which are analytic in the disk $\mathbf{D}=\{z:|z-\frac{1}{2}|\leq \frac{1}{2}\}$, including its boundary. This space is equipped with a topology of uniform convergence, and distribution on this space is any continuous linear functional. Indeed, since
\begin{eqnarray*}
\int\limits_{0}^{1}\frac{x}{1-xz}\d F_{\lambda}(x)=-\frac{\lambda}{2}G_{\lambda}(z):=\sum\limits_{L=1}^{\infty}m^{(\lambda)}_{L}z^{L-1},
\end{eqnarray*}
define a distribution $F_{\lambda}$ on the space ${\sf C}^{\omega}$ by $\langle z^{L}, F_{\lambda}\rangle=m^{\lambda}_{L}$, $L\geq 1$,  $\langle 1, F_{\lambda}\rangle=0$, and for any analytic function $B(z)\in{\sf C}^{\omega}$, $B(z)=\sum_{L=0}^{\infty}b_{L}z^{L}$, by
\begin{eqnarray*}
\langle B,F_{\lambda}\rangle=\sum\limits_{L=0}^{\infty}b_{L}\langle z^{L},F_{\lambda}\rangle.
\end{eqnarray*}
First, $\langle*,F_{\lambda}\rangle$ is certainly a linear functional and is properly defined, since the functional equation (\ref{eigen}) for $G_{\lambda}(z)$ implies that it possesses all left derivatives at $z=1$; as a consequence, the series $\sum_{L=1}^{\infty}L^{p}|m^{(\lambda)}_{L}|$ converges for any $p\in\mathbb{N}$ (see Theorem 3 for the estimates on moments $m_{L}$). Second, let $B_{n}(z)=\sum_{L=0}^{\infty}b^{(n)}_{L}z^{L}$, $n\geq 1$, converge uniformly to $B(z)$ in the circle $|z|\leq 1$. Thus, $\sup_{|z|\leq 1}|B_{n}(z)-B(z)|=r_{n}\rightarrow 0$. Then by Cauchy formula,
\begin{eqnarray*}
b^{(n)}_{L}=\frac{1}{2\pi i}\oint\limits_{|z|=1}\frac{B_{n}(z)}{z^{L+1}}\d z.
\end{eqnarray*}
This obviously imply that $|b^{(n)}_{L}-b_{L}|\leq r_{n}$, $L\geq 0$, and therefore $\langle*,F_{\lambda}\rangle$ is continuous, and hence it is a distribution.
Using the condition
$\d F_{\lambda}(x+1)=\frac{1}{2}\d F_{\lambda}(x)$, these distributions can be extended to other spaces. Summarizing, we have shown that Minkowski question mark function has an infinite sequence of ``peers" $F_{\lambda}(x)$ which are also related to continued fraction expansion, in somewhat similar manner. $F(x)$ is the only one among them being ``non-homogeneous".
\section{Linear relations among moments $M_{L}$}
In this section we clarify the nature of linear relations among moments $M_{L}$. This was mentioned in \cite{as}, but not done in explicit form. Note that (\ref{symm}) gives linear
relations among moments $m_{L}$:
$m_{L}=\sum_{s=0}^{L}\binom{L}{s}(-1)^{s}m_{s}$, $L\geq 0$. These
linear relations can be written in terms of $M_{L}$. Despite the
fact that second equality of (\ref{symm}) is general phenomena for symmetric
distribution, in conjunction with (\ref{intr}) it gives the essential
information about $F(x)$. Therefore, let us denote
\begin{eqnarray*}
{\sf
q}(x,t)=(2-e^{t})e^{xt}-(2e^{t}-1)e^{-xt}=\sum\limits_{n=1}^{\infty}{\sf
Q}_{n}(x)\frac{t^{n}}{n!}.
\end{eqnarray*}
We see that ${\sf Q}_{n}(x)$ are polynomials with integer coefficients and they are given by
\begin{eqnarray}
{\sf Q}_{n}(x)=2x^{n}-(x+1)^{n}-2(1-x)^{n}+(-x)^{n}. \label{poly}
\end{eqnarray}
The
following table gives the first few polynomials.\\

\begin{tabular}{|r | r || r | r|}
\hline
$n$ & ${\sf Q}_{n}(x)$           &$n$  & ${\sf Q}_{n}(x)$\\
\hline
$1$&  $2x-3$                  &$5$  & $2x^{5}-15x^{4}+10x^{3}-30x^{2}+5x-3$\\
$2$ & $2x-3$                  &$6$  & $6x^{5}-45x^{4}+20x^{3}-45x^{2}+6x-3$\\
$3$ & $2x^{3}-9x^{2}+3x-3$    &$7$  & $2x^{7}-21x^{6}+21x^{5}-105x^{4}+35x^{3}-63x^{2}+7x-3$\\
$4$ & $4x^{3}-18x^{2}+4x-3$   &$8$  & $8x^{7}-84x^{6}+56x^{5}-210x^{4}+56x^{3}-84x^{2}+8x-3$\\
\hline
\end{tabular}\\

Moreover, the following statement holds.
\begin{prop} Polynomials ${\sf Q}_{n}(x)$ have the following
properties:
\begin{eqnarray*}
&(i)&\quad {\sf Q}_{2n}(x)\in
L_{\mathbb{Q}}\big{(}{\sf Q}_{1}(x),{\sf Q}_{3}(x),...,{\sf Q}_{2n-1}(x)\big{)},\quad n\geq 1;\\
&(ii)&\quad \deg {\sf Q}_{2n}=2n-1,\quad \deg {\sf Q}_{2n-1}=2n-1, \quad n\geq 1;\\
&(iii)&\quad\widehat{{\sf Q}}_{2n}(x):=\frac{{\sf Q}_{2n}(x)+3}{x}
\text{ is reciprocal}:\widehat{{\sf
Q}}_{2n}(x)=x^{2n-2}\widehat{{\sf
Q}}_{2n}\Big{(}\frac{1}{x}\Big{)};\\
&(iv)&\quad \int\limits_{0}^{\infty}{\sf Q}_{n}(x)\d F(x)=0.
\end{eqnarray*}
\end{prop}
Naturally, it is property $(iv)$ which makes these polynomials very
important in the study of Minkowski's $?(x)$. Here
$L_{\mathbb{Q}}(*)$ denotes the $\mathbb{Q}-$linear space spanned
by the specified polynomials.\\

\proof{\it (i) }Let ${\sf q}_{e}(x,t)=\frac{1}{2}\big{(}{\sf q}(x,t)+{\sf
q}(x,-t)\big{)}$, and ${\sf q}_{o}(x,t)=\frac{1}{2}\big{(}{\sf q}(x,t)-{\sf
q}(x,-t)\big{)}$. Direct calculation shows that, if $e^{t}=T$, then
\begin{eqnarray*}
2{\sf
q}_{e}=e^{xt}(3-T-\frac{2}{T})+e^{-xt}(3-\frac{1}{T}-2T),\quad
2{\sf q}_{o}=e^{xt}(1-T+\frac{2}{T})-e^{-xt}(1-\frac{1}{T}+2T).
\end{eqnarray*}
This yields
\begin{eqnarray*}
\sum\limits_{n=1}^{\infty}{\sf Q}_{2n}(x)\frac{t^{2n}}{(2n)!}={\sf
q}_{e}(x,t)=\frac{T-1}{T+1}{\sf
q}_{o}(x,t)=\frac{e^{t}-1}{e^{t}+1}\sum\limits_{n=0}^{\infty}{\sf
Q}_{2n+1}(x)\frac{t^{2n+1}}{(2n+1)!}.
\end{eqnarray*}
The multiplier on the right, $\frac{e^{t}-1}{e^{t}+1}=\tanh(t/2)$, is independent of
$x$, and this obviously proves part $(i)$. Also, part
{\it (ii)} follows easily from (\ref{poly}).\\
{\it(iii)} Since $\widehat{{\sf
Q}}_{2n}(x)=\frac{1}{x}(3x^{2n}-(x+1)^{2n}-2(x-1)^{2n}+3)$, the proof is immediate.\\
{\it (iv)} In fact, (\ref{symm}) gives
$(2-e^{t})M(t)=(2e^{t}-1)M(-t)$. For real $|t|<\log 2$, we have
$M(t)=\int_{0}^{\infty}e^{xt}\d F(x)$. This implies
\begin{eqnarray*}
\int\limits_{0}^{\infty}{\sf q}(x,t)\d
F(x)=\sum\limits_{n=0}^{\infty}\frac{t^{n}}{n!}\int\limits_{0}^{\infty}{\sf
Q}_{n}(x)\d F(x)\equiv0,\quad\text{for }|t|<\log 2,
\end{eqnarray*}
and this completes the proof. $\blacksquare$\\
Consequently, there exist linear relations among moments $M_{L}$. Thus, for
example, part $(iv)$ (in case $n=1$ and $n=3$) implies $2M_{1}-3=0$ and
$2M_{3}-9M_{2}+3M_{1}=3$ respectively. The exact values of $M_{L}$ belong to the class of
constants, which can be thought as emerging from arithmetic-geometric chaos.
This resembles the situation concerning polynomial relations among various
periods. We will not present the definition of a period, which can be found in
\cite{konzag}. In particular, the authors conjecture (and there is no support
for possibility that it can be proved wrong) that \it ``if a period has two
integral representations, then one can pass from one formula to another using
only additivity, change of variables, and Newton-Leibniz formula, in which all
functions and domains of integration are algebraic with coefficients in
$\overline{\mathbb{Q}}$". \rm Thus, for example, the conjecture predicts the
possibility to prove directly that $ \iint_{\frac{x^{2}}{4}+3y^{2}\leq 1}\d x\d
y=\int_{-1}^{1} \frac{\d x}{\sqrt[3]{(1-x)(1+x)^{2}}}, $ without knowing that
they both are equal to $\frac{2\pi}{\sqrt{3}}$, and this indeed can be
done. Similarly, returning to the topic of this paper, we believe that
any finite $\mathbb{Q}-$linear relation among constants $M_{L}$ can be proved
simply by applying the functional equation of $F(x)$, by means of integration by
parts and change of variables. The last proposition supports this claim. In other words, we believe that there cannot be any other miraculous
coincidences regarding the values
of $M_{L}$. More precisely, \\

\textbf{Conjecture. }\it Suppose, $r_{k}\in\mathbb{Q}$, $0\leq
k\leq L$, are rational numbers such that
$\sum_{k=0}^{L}r_{k}M_{k}=0$. Let
$\ell=\big{[}\frac{L-1}{2}\big{]}$. Then
\begin{eqnarray*}
\sum\limits_{k=0}^{L}r_{k}x^{k}\in L_{\mathbb{Q}}\big{(}{\sf
Q}_{1}(x),{\sf Q}_{3}(x),...,{\sf Q}_{2\ell+1}(x)\big{)}.
\end{eqnarray*}
\rm This conjecture, if true, should be difficult to prove. It would imply,
for example, that $M_{L}$ for $L\geq 2$ are irrational. On the other hand, this
conjecture seems to be much more natural and approachable, compared to similar
conjectures regarding arithmetic nature of constants emerging from geometric
chaos, e.g. spectral values $s$ for Maass wave forms (say, for ${\sf
PSL}_{2}(\mathbb{Z})$), or those coming from arithmetic chaos, like non-trivial
zeros of Riemann's $\zeta(s)$. We cannot give any other evidence, save the last
proposition, to support this conjecture.
\section{Estimate for the moments $m_{L}$}
This chapter deals with the asymptotic estimate for the moments $m_{L}$. This result was not obtained before, and in view of the expression (\ref{mom}), it is of certain number-theoretic interest. This result should be compared with the asymptotic formula for $M_{L}$, see (\ref{sum}) and the expression just below it. {\it A priori}, as it is implied by the fact that the radius of convergence of $G(z)$ at $z=0$ is 1, and by (\ref{est}), for every $\varepsilon>0$ and $p>1$, one has $\frac{1}{L^{p}}\gg m_{L}\gg(1-\varepsilon)^{L}$, as $L\rightarrow\infty$. More precisely, we have
\begin{thm}Let $C=e^{-2\sqrt{\log2}}=0.18917...$. Then the following estimate holds, as $L\rightarrow\infty$:
\begin{eqnarray*}
C^{\sqrt{L}}\ll m_{L}\ll L^{1/4}C^{\sqrt{L}}.
\end{eqnarray*}
Both implied constants are absolute.
\end{thm}
\proof Fix $J\in\mathbb{N}$, and choose an increasing sequence of positive real numbers $\mu_{j}<1$, $1\leq j\leq J$. We will soon specify $\mu_{j}$ in such a way that $\mu_{j}\rightarrow0$ uniformly as $L\rightarrow\infty$. An estimate for $m_{L}$ is obtained via the defining integral (recall that $F(x)+F(1/x)=1$):
\begin{eqnarray*}
m_{L}=(\int\limits_{0}^{\mu_{1}}+\sum\limits_{j=1}^{J-1}\int\limits_{\mu_{j}}^{\mu_{j+1}}+\int\limits_{\mu_{J}}^{\infty})
\Big{(}\frac{1}{x+1}\Big{)}^{L}\d F(x)<\\
F(\mu_{1})+\sum\limits_{j=1}^{J-1}\Big{(}\frac{1}{\mu_{j}+1}\Big{)}^{L}F(\mu_{j+1})+\Big{(}\frac{1}{\mu_{J}+1}\Big{)}^{L}.\\
\end{eqnarray*}
Indeed, in the first integral, the integrand is bounded by $1$. In the middle integrals, we choose the largest value of integrand, and change bounds of integration to $[0,\mu_{j+1}]$. The same is done with the last integral, with bounds changed to $[0,\infty)$.  Now choose $\mu_{j}=\frac{1}{c_{j}\sqrt{L}}$ for some decreasing sequence of constants $c_{j}$. Functional equation for $F(x)$ implies $F(x+n)=1-2^{-n}+2^{-n}F(x)$, $x\geq 0$. Thus, $1-F(x)\asymp 2^{-x}$, as $x\rightarrow\infty$ (the implied constants being $\min$ and $\max$ of the function $\Psi(x)$; see Figure 3 and Section 8). Using the identity $F(x)+F(1/x)=1$, we therefore obtain
\begin{eqnarray}
m_{L}\ll2^{-c_{1}\sqrt{L}}+
\sum\limits_{j=1}^{J-1}\Big{(}\frac{1}{\frac{1}{c_{j}\sqrt{L}}+1}\Big{)}^{L}2^{-c_{j+1}\sqrt{L}}+
\Big{(}\frac{1}{\frac{1}{c_{J}\sqrt{L}}+1}\Big{)}^{L}\ll\nonumber\\
e^{-\sqrt{L}{c_{1}\log2}}+\sum\limits_{j=1}^{J-1}e^{-\sqrt{L}(\frac{1}{c_{j}}+c_{j+1}\log 2)}+e^{-\sqrt{L}\frac{1}{c_{J}}}.\label{estim}
\end{eqnarray}
Here we need an elementary lemma.
\begin{lem} For given $J\in\mathbb{N}$, there exists a unique sequence of positive real numbers $c^{*}_{1},...,c^{*}_{J}$, such that
\begin{eqnarray*}
c^{*}_{1}=\frac{1}{c^{*}_{1}}+c^{*}_{2}=\frac{1}{c^{*}_{2}}+c^{*}_{3}=...=\frac{1}{c^{*}_{J-1}}+c^{*}_{J}=\frac{1}{c^{*}_{J}}.
\end{eqnarray*}
Moreover, this sequence $\{c^{*}_{j},1\leq j\leq J\}$ is decreasing, and it is given by
\begin{eqnarray*}
c_{j}^{*}=\frac{\sin\frac{(j+1)\pi}{J+2}}{\sin\frac{j\pi}{J+2}},\quad j=1,2,...,J\Rightarrow c_{1}^{*}=2\cos \frac{\pi}{J+2}.
\end{eqnarray*}
\end{lem}
{\it Proof. }Indeed, we see that $c^{*}_{1}=x$ determines the sequence $c^{*}_{j}$ uniquely. First, $c_{2}=x-\frac{1}{x}=\frac{x^{2}-1}{x}$. Let $F_{1}(x)=x$, $F_{2}(x)=x^{2}-1$. Suppose we have shown that $c_{j}=\frac{F_{j}(x)}{F_{j-1}(x)}$ for certain sequence of polynomials. Then from the above equations one obtains
\begin{eqnarray*}
c_{j+1}=c_{1}-\frac{F_{j-1}(x)}{F_{j}(x)}=\frac{xF_{j}(x)-F_{j-1}(x)}{F_{j}(x)}.
\end{eqnarray*}
Thus, using induction we see that $c_{j}=\frac{F_{j}(x)}{F_{j-1}(x)}$, where polynomials $F_{j}(x)$ are given by the initial values $F_{0}(x)=1$, $F_{1}(x)=x$ and then for $j\geq 1$ recurrently by $F_{j+1}(x)=xF_{j}(x)-F_{j-1}(x)$. This shows that $F_{j}(2x)=U_{j}(x)$, where $U(x)$ stand for the classical Chebyshev $U-$polynomials, given by
\begin{eqnarray*}
U_{j}(\cos \theta)=\frac{\sin(j+1)\theta}{\sin \theta}.
\end{eqnarray*}
The last equation $c_{1}^{*}=\frac{1}{c^{*}_{J}}$ implies
$F_{J+1}(x)=0$. Thus, $U_{J+1}(x/2)=0$, and all possible values of $c^{*}_{1}$ are given by $c_{1}^{*}=x=2\cos \frac{k\pi}{J+2}$, $k=1,2,...,J+1$. Thus,
\begin{eqnarray*}
c_{j}^{*}=\frac{F_{j}(x)}{F_{j-1}(x)}=\frac{U_{j}(x/2)}{U_{j-1}(x/2)}=\frac{\sin\frac{k(j+1)\pi}{J+2}}{\sin\frac{kj\pi}{J+2}}.
\end{eqnarray*}
Since our concern is only positive solutions, this gives the last statement of lemma. Finally, monotonicity is easily verifiable. Indeed, system of equations imply $c^{*}_{2}<c^{*}_{1}$, and then we act by induction. $\square$\\
Thus, $c_{1}^{*}>2-\frac{b}{J^{2}}$ for some constant $b>0$.
Returning to the proof of the Theorem, for given $J$, let $c^{*}_{j}$ be the sequence in lemma, and let $c^{*}_{i}=c_{i}\sqrt{\log2}$. Thus,
\begin{eqnarray*}
c_{1}\log2=\frac{1}{c_{1}}+c_{2}\log2=\frac{1}{c_{2}}+c_{3}\log2=...=\frac{1}{c_{J-1}}+c_{J}\log2=\frac{1}{c_{J}}.
\end{eqnarray*}
Choosing exactly this sequence for the estimate (\ref{estim}), and using the bound for $c^{*}_{1}$, we get:
\begin{eqnarray*}
m_{L}\ll(J+1)e^{-\sqrt{L}c_{1}\log 2}<(J+1)C^{\sqrt{L}}e^{\frac{b\sqrt{\log2}}{J^{2}}\sqrt{L}}.
\end{eqnarray*}
Finally, the choice $J=[L^{1/4}]$ establishes the upper bound.\\
The lower estimate is immediate. In fact, let $\mu=\frac{1}{c\sqrt{L}}$. Then
\begin{eqnarray*}
m_{L}>\int\limits_{0}^{\mu}\Big{(}\frac{1}{x+1}\Big{)}^{L}\d F(x)>\Big{(}\frac{1}{\mu+1}\Big{)}^{L}F(\mu)\gg2^{-c\sqrt{L}}\cdot e^{-\sqrt{L}\frac{1}{c}}
\end{eqnarray*}
The choice $c=\log^{-1/2}2$ gives the desired bound.$\blacksquare$\\
The constants in Theorem can also be calculated without great effort, but this is astray of the main topic of the paper.\\
It should be noted that, if we start directly from the second definition
(\ref{moments}) of $m_{L}$, then in the course of the proof of Theorem 3, we use both equalities $F(x)+F(1/x)=1$ and $2F(\frac{x}{x+1})=F(x)$. Since these two determine $F(x)$ uniquely, generally speaking, our estimate for $m_{L}$ is characteristic only to $F(x)$. A direct inspection of the proof also reveals that the true asymptotic ``action" in the second definition (\ref{moments}) of $m_{L}$ takes place in the neighborhood of $1$. This, obviously, is a general fact for probabilistic distributions with proper support on the interval $[0,1]$. Additionally, calculations show that the sequence $m_{L}/(L^{1/4}C^{\sqrt{L}})$ is monotonically increasing. If this is indeed the case, there exists a limit $A=\lim_{L\rightarrow\infty}\frac{m_{L}}{L^{1/4}C^{\sqrt{L}}}$.
Such a constant $A$ might be of definite interest in metric number theory. This is the main topic of the paper \cite{ga4}.\\

As a final remark, we note that the result of Theorem must be considered in conjunction with the linear relations $m_{L}=\sum_{s=0}^{L}\binom{L}{s}(-1)^{s}m_{s}$, $L\geq 0$, and the natural inequalities, imposed by the fact that $m_{L}$ is a sequence of moments of probabilistic distribution with support on the interval $[0,1]$. We thus have Hausdorff conditions, which state that for all non-negative integers $m$ and $n$, one has
\begin{eqnarray*}
2\int\limits_{0}^{1}x^{n}(1-x)^{m}\d F(x)=\sum\limits_{i=0}^{m}\binom{m}{i}(-1)^{i}m_{i+n}>0.
\end{eqnarray*}
This is, of course, the consequence of monotonicity of $F(x)$.
\section{Exact sequence}
In this section, we prove the exactness of a sequence of continuous linear maps, intricately related to Minkowski's question mark function $F(x)$. Let ${\sf C}[0,1]$ denote the space of continuous, complex-valued functions on the interval $[0,1]$ with supremum norm. For $f\in{\sf C}[0,1]$, one has the identity
\begin{eqnarray}
\int\limits_{0}^{1}f(x)\d F(x)=\sum\limits_{n=1}^{\infty}
\int\limits_{0}^{1}f\Big{(}\frac{1}{x+n}\Big{)}2^{-n}\d F(x),\label{summa}
\end{eqnarray}
Indeed, using functional equation (\ref{distr}), we have
\begin{eqnarray*}
\int\limits_{0}^{1}f(x)\d F(x)=\int\limits_{1}^{\infty}f\Big{(}\frac{1}{x}\Big{)}\d F(x)=
\sum\limits_{n=1}^{\infty}\int\limits_{0}^{1}f\Big{(}\frac{1}{x+n}\Big{)}\d F(x+n),
\end{eqnarray*}
which is exactly (\ref{summa}). Let ${\sf C}^{\omega}$ denote, as before, the space of analytic functions in the disk $\mathbf{D}=|z-\frac{1}{2}|\leq\frac{1}{2}$, including its boundary. We equip this space with the topology of uniform convergence (as a matter of fact, we have a wider choice of spaces; this one is chosen as an important example). Now, consider a continuous functional on ${\sf C}^{\omega}$ given by $T(f)=\int_{0}^{1}f(x)\d F(x)$, and a continuous non-compact linear operator
$
[\mathcal{L}f](x)=f(x)-\sum_{n=1}^{\infty}f\big{(}\frac{1}{x+n}\big{)}2^{-n}$.
Finally, let $i$ stand for the natural inclusion $i:\mathbb{C}\rightarrow {\sf C}^{\omega}$.
\begin{thm}
The following sequence of maps is exact:
\begin{eqnarray*}
0\rightarrow\mathbb{C}\mathop{\rightarrow}^{i} \mathop{{\sf C}^{\omega}}_{(*)}\mathop{\rightarrow}^{\mathcal{L}} \mathop{{\sf C}^{\omega}}_{(**)}\mathop{\rightarrow}^{T}\mathbb{C}\rightarrow0.
\end{eqnarray*}
\end{thm}
\proof First, $i$ is obviously a monomorphism. Let $f\in\text{Ker}(\mathcal{L})$. This means that
$f(x)=\sum_{n=1}^{\infty}f\big{(}\frac{1}{x+n}\big{)}2^{-n}$. Let $x_{0}\in[0,1]$ be such that $|f(x_{0})|=\mathop{\sup}\limits_{x\in[0,1]}|f(x)|$. Since $\sum_{n=1}^{\infty}2^{-n}=1$, this yields $f(\frac{1}{x_{0}+n})=f(x_{0})$ for $n\in\mathbb{N}$. By induction, $f([0,n_{1},n_{2},...,n_{I}+x_{0}])=f(x_{0})$ for all $I\in\mathbb{N}$, and all $n_{i}\in\mathbb{N}$, $1\leq i\leq I$; here $[*]$ stands for the (regular) continued fraction. Since this set is everywhere dense in $[0,1]$ and $f$ is continuous, this forces $f(x)\equiv\text{const}$ for $x\in[0,1]$. Due to the analytic continuation, this is valid for $x\in\mathbf{D}$ as well. Hence, we have the exactness at term $(*)$.\\
Next, $T$ is obviously an epimorphism. Further, identity (\ref{summa}) implies that
$\text{Im}(\mathcal{L})\subset\text{Ker}(T)$. The task is to show that indeed we have an equality. At this stage, we need the following lemma. Denote $[\mathcal{S}f](x)=\sum_{n=1}^{\infty}f\big{(}\frac{1}{x+n}\big{)}2^{-n}$.
\begin{lem} Let $f\in{\sf C}^{\omega}$.
Then $[\mathcal{S}^{n}f](x)=2T(f)+O(\gamma^{-2n})$ for $x\in\mathbf{D}$.\\

Here $T(f)$ stands for the constant function, $\gamma=\frac{1+\sqrt{5}}{2}$ is the golden section, and the bound implied by $O$ is uniform for $x\in\mathbf{D}$.
\end{lem}
\proof In fact, lemma is true for any function with continuous derivative. Let $x\in\mathbf{D}$. We have
\begin{eqnarray*}
[\mathcal{S}^{r}f](x)=\sum\limits_{n_{1},n_{2},...,n_{r}=1}^{\infty}2^{-(n_{1}+n_{2}+...+n_{r})}f([0,n_{1},n_{2},...,n_{r}+x]).
\end{eqnarray*}
The direct inspection of this expression and (\ref{min}) shows that this is exactly twice the Riemann sum for the integral $\int\limits_{0}^{1}f(x)\d F(x)$, corresponding to the division of unit interval into intervals with endpoints being $[0,n_{1},n_{2},...,n_{r}]$, $n_{i}\in\mathbb{N}$. From the basic properties of M\"{o}bius transformations we inherit that the set $[0,n_{1},n_{2},...,n_{r}+x]$ for $x\in\mathbf{D}$ is a circle $\mathbf{D}_{r}$ whose diagonal is one of these intervals, say $I_{r}$. For fixed $r$, the largest of these intervals has endpoints $\frac{F_{r-1}}{F_{r}}$ and $\frac{F_{r}}{F_{r+1}}$, where $F_{r}$ stands for the usual Fibonacci sequence. Thus, its length is $\frac{1}{F_{r}F_{r+1}}\sim c\gamma^{-2r}$. Let $x_{0},x_{1}\in \mathbf{B}_{r}$, and $\sup_{x\in\mathbf{D}}|f'(x)|=A$. We have
\begin{eqnarray*}
\sup_{x_{0},x_{1}\in\mathbf{D}_{r}}|f(x_{0})-f(x_{1})|\leq Ac\gamma^{-2r}.
\end{eqnarray*}
Thus, the Riemann sum deviates from the Riemann integral no more than
\begin{eqnarray*}
|[\mathcal{S}^{r}f](x)-2T(f)|\leq Ac\gamma^{-2r}\sum\limits_{n_{1},n_{2},...,n_{r}=1}^{\infty}2^{-(n_{1}+n_{2}+...+n_{r})}=Ac\gamma^{-2r}. \end{eqnarray*}
This proves lemma. $\square$\\
Thus, let $f\in\text{Ker}(T)$. All we need is to show that the equation $f=g-\mathcal{S}g$ has a solution $g\in{\sf C}^{\omega}$. Indeed, let
$g=f+\sum_{n=1}^{\infty}\mathcal{S}^{n}f$. By the above lemma, $\|\mathcal{S}^{n}f\|=O(\gamma^{-2n})$. Thus, the series defining $g$ converges uniformly and hence $g$ is an analytic function. Finally, $g-\mathcal{S}g=f$; this shows that $\text{Ker}(T)\subset\text{Im}(\mathcal{L})$ and exactness at term $(**)$ is proved.$\blacksquare$\\

The eigen-functions of $\mathcal{S}$ acting on the space ${\sf C}^{\omega}$ are given by
$G^{\star}(-x)=\int_{0}^{-x}G_{\lambda}(z)\d z+\int_{-1}^{0}G_{\lambda}(z)\d z$ (see equations (\ref{summa2}) and (\ref{summa3}) in the next section). Thus, the problem of convergence of $\mathcal{S}^{n}f$ is completely analogous to the problem of convergence for the iterates of Gauss-Kuzmin-Wirsing operator. Let us remind that if $f\in{\sf C}[0,1]$, it is given by
$
[\mathcal{W}f](x)=\sum_{n=1}^{\infty}\frac{1}{(x+n)^{2}}f\big{(}
\frac{1}{x+n}\big{)}.$
Dominant eigen-value $1$ correspond to an eigen-function $\frac{1}{1+x}$. As it was proved by Kuzmin, provided that $f(x)$ has a continuous derivative, there exists $c>0$, such that
\begin{eqnarray*}
[\mathcal{W}^{n}f](x)=\frac{A}{1+x}+O(e^{-c\sqrt{n}}),\text{ as }n\rightarrow\infty;\quad A=\frac{1}{\log 2}\int\limits_{0}^{1}f(x)\d x.
\end{eqnarray*}
The proof can be found in \cite{khin}. Note that this was already conjectured by Gauss, but he did not give the proof nor for the main neither for the error term. For the most important case, when $f(x)=1$, L\'{e}vy established the error term of the form $O(C^{n})$ for $C=0.7$. Finally, Wirsing \cite{wirsing} gave the exact result in terms of eigen-functions of $\mathcal{W}$, establishing the error term of the form $(-c)^{n}\Psi(x)+O(x(1-x)\mu^{n})$, where $c=-0.303663...$ is the sub-dominant eigen-value (Gauss-Kuzmin-Wirsing constant), $\Psi(x)$ is a corresponding eigen-function, and $\mu<|c|$. Returning to our case, we have completely analogous situation: operator $\mathcal{W}$ is replaced by $\mathcal{S}$, and the measure $\d x$ is replaced by $\d F(x)$. The leading eigen-value $1$ corresponds to the constant function. However strange, Wirsing did not notice that eigen-values of $\mathcal{W}$ are in fact eigen-values of certain Hilbert-Schmidt operator. This was later clarified by Bobenko \cite{bobenko}. Recently, Gauss-Kuzmin-L\'{e}vy theorem was generalized by Manin and Marcolli in \cite{manin}. The paper is very rich in ideas and results; in particular, in sheds a new light on the theorem just mentioned.\\

Concerning spaces for which Theorem 4 holds, we can investigate space the ${\sf C}[0,1]$ as well. However, if $f\in{\sf C}[0,1]$ and $f\in\text{Ker}(T)$, the significant difficulty arises in proving uniform convergence of the series $\sum_{n=0}^{\infty}\mathcal{S}^{n}f$. Moreover, operator $\mathcal{S}$, acting on the space ${\sf C}[0,1]$, has additional point spectrum apart from $\lambda$. Indeed, let $P_{n}(y)=y^{n}+\sum\limits_{i=0}^{n-1}a_{i}y^{i}$ be a polynomial of degree $n$, which satisfies yet another variation of three term functional equation
\begin{eqnarray*}
2P_{n}(1-2y)-P_{n}(1-y)=\frac{1}{\delta_{n}}P(y)
\end{eqnarray*}
for certain $\delta_{n}$. The comparison of leading terms shows that $\delta_{n}=\frac{(-1)^{n}}{2^{n+1}-1}$, and that indeed for this $\delta_{n}$ there exists a unique polynomial, since each coefficient $a_{j}$ can be uniquely determined with the knowledge of coefficients $a_{i}$ for $i>j$. Thus, \\

\begin{tabular}{l l}
$P_{1}(y)=y-\frac{1}{4}$,             & $P_{2}(y)=y^{2}-\frac{3}{5}y+\frac{1}{15}$,\\
$P_{3}(y)=y^{3}-\frac{21}{22}y^{2}+\frac{3}{11}y-\frac{7}{352}$, & $P_{4}(y)=y^{4}-\frac{30}{23}y^{3}+\frac{14}{23}y^{2}-\frac{45}{391}y +\frac{37}{5865}$.
\end{tabular}\\

 The equation for $P_{n}$ implies that
\begin{eqnarray}
\delta_{n}P_{n}(y)=\sum\limits_{n=1}^{\infty}\frac{1}{2^{n}}P_{n}\Big{(}\frac{1-y}{2^{n}}\Big{)}.\label{pol}
\end{eqnarray}
Then we have
\begin{prop}
The function $P_{n}(F(x))$ is the eigen-function of $\mathcal{S}$, acting on the space ${\sf C}[0,1]$, and eigen-value $\frac{(-1)^{n}}{2^{n+1}-1}$ belongs to the point spectrum of $\mathcal{S}$.
\end{prop}
\proof
Indeed,
\begin{eqnarray*}
[\mathcal{S}(P_{n}\circ F)](x)=\sum\limits_{n=1}^{\infty}\frac{1}{2^{n}}P_{n}\circ F\Big{(}\frac{1}{x+n}\Big{)}\mathop{=}^{(\ref{distr})}
\frac{1}{2^{n}}P_{n}\Big{(}1-F(x+n)\Big{)}\mathop{=}^{(\ref{distr})}\\
\sum\limits_{n=1}^{\infty}\frac{1}{2^{n}}P_{n}\Big{(}2^{-n}-2^{-n}F(x)\Big{)}\mathop{=}^{(\ref{pol})}
\delta_{n}P_{n}(F(x)).\blacksquare
\end{eqnarray*}
Thus, operator $\mathcal{S}$ behaves differently in spaces ${\sf C}[0,1]$ and ${\sf C}^{\omega}$. We postpone the analysis of this operator in various spaces for the consequent paper.
\section{Integrals, involving $F(x)$}
In this section we calculate certain integrals. Only rarely it is possible to express an integral involving $F(x)$ in closed form. In fact, all results we possess come from the
identity $M_{1}=\frac{3}{2}$, and any iteration of identities similar to (\ref{summa}). The
following theorem adds identities of quite a different sort.
\begin{thm} Let $G_{\lambda}(z)$ be any function, which satisfies
the hypotheses of Theorem 1. Then
\begin{eqnarray*}
&(i)&\frac{\lambda}{\lambda+1}\int\limits_{0}^{1}G_{\lambda}(-x)\d x=
\int\limits_{0}^{1}G_{\lambda}(-x)F(x)\d x;\\
&(ii)&-\int\limits_{0}^{1}\log x\d F(x)=2\int\limits_{0}^{1}\log(1+x)\d F(x)=\int\limits_{0}^{1}G(-x)\d x;\\
&(iii)& \int\limits_{0}^{1}G(-x)(1+x^{2})\d F(x)=\frac{1}{4};\\
&(iv)& \int\limits_{0}^{1}G_{\lambda}(-x)\Big{(}1-\frac{x^{2}}{\lambda}\Big{)}\d F(x)=0.
\end{eqnarray*}
\end{thm}
\proof We first prove identity $(i)$. By (\ref{eigen}), for every integer $n\geq 1$, we have
\begin{eqnarray*}
2G_{\lambda}(-z-n+1)-G_{\lambda}(-z-n)=\frac{1}{\lambda (z+n)^{2}}G_{\lambda}\Big{(}-\frac{1}{z+n}\Big{)}.
\end{eqnarray*}
Divide this by $2^{n}$ and, sum over $n\geq 1$. By Theorem 2, the sum on the left is absolutely convergent. Thus,
\begin{eqnarray*}
G_{\lambda}(-z)=\sum\limits_{n=1}^{\infty}\frac{1}{\lambda2^{n} (z+n)^{2}}G_{\lambda}\Big{(}-\frac{1}{z+n}\Big{)}
\end{eqnarray*}
Let $G^{\star}_{\lambda}(x)=\int_{0}^{x}G_{\lambda}(z)\d z$.
In terms of $G^{\star}_{\lambda}(x)$, the last identity reads as
\begin{eqnarray}
-G^{\star}_{\lambda}(-x)=\sum\limits_{n=1}^{\infty}\frac{1}{\lambda2^{n}}
G^{\star}_{\lambda}\Big{(}-\frac{1}{x+n}\Big{)}-\sum\limits_{n=1}^{\infty}\frac{1}{\lambda2^{n}}
G^{\star}_{\lambda}\Big{(}-\frac{1}{n}\Big{)}.\label{summa2}
\end{eqnarray}
In particular, setting $x=1$, one obtains
\begin{eqnarray}
\sum\limits_{n=1}^{\infty}\frac{1}{\lambda2^{n}}
G^{\star}_{\lambda}\Big{(}-\frac{1}{n}\Big{)}=(\frac{1}{\lambda}-1)G^{\star}_{\lambda}(-1).\label{summa3}
\end{eqnarray}
Now we are able to calculate the following integral (we use integration by parts in Stieltjes integral twice).
\begin{eqnarray*}
\int\limits_{0}^{1}G_{\lambda}(-x)F(x)\d x=-\int\limits_{0}^{1}\frac{\d}{\d x}G^{\star}_{\lambda}(-x)F(x)\d x=-\frac{1}{2}G^{\star}_{\lambda}(-1)+
\int\limits_{0}^{1}G^{\star}_{\lambda}(-x)\d F(x)\mathop{=}^{(\ref{summa2})}\\
-\frac{1}{2}G^{\star}_{\lambda}(-1)+\frac{1}{2}\sum\limits_{n=1}^{\infty}\frac{1}{\lambda2^{n}}
G^{\star}_{\lambda}\Big{(}-\frac{1}{n}\Big{)}-\frac{1}{\lambda}\sum\limits_{n=1}^{\infty}
\int\limits_{0}^{1}G^{\star}_{\lambda}\Big{(}-\frac{1}{x+n}\Big{)}2^{-n}\d F(x)\mathop{=}^{(\ref{summa}),(\ref{summa3})}\\
-\frac{1}{2}G^{\star}(-1)+\frac{1}{2}(\frac{1}{\lambda}-1)G^{\star}_{\lambda}(-1)-
\frac{1}{\lambda}\int\limits_{0}^{1}G^{\star}_{\lambda}(-x)\d F(x)=-G^{\star}(-1)-\frac{1}{\lambda}\int\limits_{0}^{1}G_{\lambda}(-x)F(x)\d x.
\end{eqnarray*}
Thus, the same integral is on the both sides, and
\begin{eqnarray*}
\int\limits_{0}^{1}G_{\lambda}(-x)F(x)\d x=-\frac{\lambda}{\lambda+1}G^{\star}_{\lambda}(-1).
\end{eqnarray*}
This establishes the statement $(i)$. \\
Now we proceed with second identity. \ignore{One can try to imitate the last trick in case $G_{\lambda}(x)$ is replaced by $G(x)$. Unfortunately, since in this case we are dealing with "eigen-value" $-1$, the integrals involving $F(x)$ do cancel (as it could be anticipated from the first identity and the presence of multiplier $\lambda$+1). Since $(\ref{sim})$ is not homogeneous, one could  still expect that the method leads to an interesting identity. Nevertheless, everything it produces is the second identity of the theorem, which can be proved directly.
The proof of the second part mimics the previous one, but since functional equation (\ref{sim}) is not a homogeneous one, we should be careful, and thus we indulge in repeating the same (slightly modified) argument again. Using (\ref{sim}), one obtains
\begin{eqnarray*}
G(-z)=-\sum\limits_{n=1}^{\infty}\frac{v'(z+n)}{2^{n}}-\sum\limits_{n=1}^{\infty}
\frac{1}{2^{n}(z+n)^{2}}G\Big{(}-\frac{1}{z+n}\Big{)}.
\end{eqnarray*}
Let $G^{\star}(x)=\int_{0}^{x}G(z)\d z$. Along similar line,
\begin{eqnarray}
-G^{\star}(-x)=-\sum\limits_{n=1}^{\infty}\frac{v(x+n)}{2^{n}}+\sum\limits_{n=1}^{\infty}\frac{v(n)}{2^{n}}
-
\sum\limits_{n=1}^{\infty}\frac{1}{2^{n}}
G^{\star}\Big{(}-\frac{1}{x+n}\Big{)}+\sum\limits_{n=1}^{\infty}\frac{1}{2^{n}}
G^{\star}\Big{(}-\frac{1}{n}\Big{)},\label{summa4}
\end{eqnarray}
Where $C=\sum_{n=1}^{\infty}2^{-n}\log n$. Specifying $x=-1$, one obtains
\begin{eqnarray}
\sum\limits_{n=1}^{\infty}\frac{1}{2^{n}}
G^{\star}\Big{(}-\frac{1}{n}\Big{)}=\frac{v(1)}{2}-\sum\limits_{n=2}^{\infty}\frac{v(n)}{2^{n}}+2G^{\star}(-1).\label{summa5}
\end{eqnarray}
In the similar fashion, using integration by parts and previous expressions, one obtains:
\begin{eqnarray*}
\int\limits_{0}^{1}G(-x)F(x)\d x=-\int\limits_{0}^{1}\frac{\d}{\d x}G^{\star}(-x)F(x)\d x=-\frac{1}{2}G^{\star}(-1)+
\int\limits_{0}^{1}G^{\star}(-x)\d F(x)\mathop{=}^{(\ref{summa4})}\\
-\frac{1}{2}G^{\star}_{\lambda}(-1)-\frac{1}{2}\sum\limits_{n=1}^{\infty}\frac{1}{2^{n}}
G^{\star}\Big{(}-\frac{1}{n}\Big{)}-\frac{C}{2}+\\
\sum\limits_{n=1}^{\infty}\int\limits_{0}^{1}\frac{v(x+n)}{2^{n}}\d F(x)+\sum\limits_{n=1}^{\infty}
\int\limits_{0}^{1}G^{\star}\Big{(}-\frac{1}{x+n}\Big{)}2^{-n}\d F(x)\mathop{=}^{(\ref{summa}),(\ref{summa5})}\\
-\frac{3}{2}G^{\star}(-1)-\frac{v(1)}{2}+\int\limits_{0}^{1}v(\frac{1}{x})\d F(x)+
\int\limits_{0}^{1}G^{\star}(-x)\d F(x)=\\
-G^{\star}(-1)-\frac{v(1)}{2}+\int\limits_{0}^{1}v (\frac{1}{x})\d F(x) +\int\limits_{0}^{1}G(-x)F(x)\d x.
\end{eqnarray*}
Thus, integrals on both sides cancel (this could be anticipated minding the identity in the first part and the fact that now we are dealing with the case of "eigen-value" $\lambda=-1$). Thus, this gives the value of integral $\int_{0}^{1}G(-x)\d x$.} Integral (\ref{gen}) and Fubini theorem imply
\begin{eqnarray*}
\int\limits_{0}^{1}G(-z)\d z=2\int\limits_{0}^{1}\int\limits_{0}^{1}\frac{x}{1+xz}\d z\d F(x)=2\int\limits_{0}^{1}\log(1+x)\d F(x).
\end{eqnarray*}
Lastly, we apply (\ref{summa}) twice to obtain the needed equality. Indeed,
\begin{eqnarray*}
I=\int\limits_{0}^{1}\log(1+x)\d F(x)\mathop{=}^{(\ref{summa})}\sum\limits_{n=1}^{\infty}\frac{1}{2^{n}}
\int\limits_{0}^{1}\log\Big{(}1+\frac{1}{x+n}\Big{)}\d F(x)=\\
\sum\limits_{n=1}^{\infty}\frac{1}{2^{n}}\int\limits_{0}^{1}\log(x+n)\d F(x)-I\mathop{=}^{(\ref{summa})}-\int\limits_{0}^{1}\log x\d F(x)-I.
\end{eqnarray*}
This finishes the proof of $(ii)$.\\
In proving $(iii)$, we can be more concise, since the pattern of the proof goes along the same line. One has
\begin{eqnarray*}
G(-z)=-\sum\limits_{n=1}^{\infty}\frac{1}{2^{n}(z+n)^{2}}G\Big{(}-\frac{1}{z+n}\Big{)}+
\sum\limits_{n=1}^{\infty}\frac{1}{2^{n}(z+n)}.
\end{eqnarray*}
Thus,
\begin{eqnarray*}
\int\limits_{0}^{1}G(-x)\d F(x)=-\sum\limits_{n=1}^{\infty}\int\limits_{0}^{1}\frac{1}{2^{n}(x+n)^{2}}G\Big{(}-\frac{1}{x+n}\Big{)}\d F(x)+\\
\sum\limits_{n=1}^{\infty}\int\limits_{0}^{1}\frac{1}{2^{n}(x+n)}\d F(x)\mathop{=}^{(\ref{summa})}
-\int\limits_{0}^{1}x^{2}G(-x)\d F(x)+\int\limits_{0}^{1}x\d F(x).
\end{eqnarray*}
Since $\int_{0}^{1}x\d F(x)=\frac{m_{1}}{2}=\frac{1}{4}$, this finishes the proof of $(iii)$. Part $(iv)$ is completely analogous. $\blacksquare$\\
Part $(iii)$, unfortunately, gives no new information about $m_{L}$. Indeed, the identity can be rewritten as
\begin{eqnarray*}
\sum\limits_{L=1}^{\infty}m_{L}(-1)^{L-1}(m_{L-1}+m_{L+1})=\frac{1}{2},
\end{eqnarray*}
which, after regrouping, turns into the identity $m_{0}m_{1}=\frac{1}{2}$.\\
Concerning part $(iv)$, and taking into account Theorem 4, one could expect that in fact $\text{Ker}(T)$ is equal to the closure of vector space spanned by functions $G_{\lambda}(-x)(1-\frac{x^{2}}{\lambda})$. If this is the case, then these functions, along with $G(z)(1+x^{2})$, produce a Schauder basis for ${\sf C}^{\omega}$. Thus, if $x^{L}=\sum_{\lambda}a^{(\lambda)}_{L}G_{\lambda}(-x)(1-\frac{x^{2}}{\lambda})$, then $a^{(-1)}_{L}=2m_{L}$.
We hope to return to this point in a future paper.\\
Concerning $(i)$, note that the values of both integrals depend on the
normalization of $G_{\lambda}$, since it is an eigen-function.
Replacing $G_{\lambda}(z)$by $cG_{\lambda}(z)$ for some
$c\in\mathbb{R}$, we deduce that the left integral is equal to
$1$ or $0$. Then $(i)$ states that
$\int\limits_{0}^{1}F(x)G_{\lambda}(-x)\d
x=\frac{\lambda}{\lambda+1}\text{ or } 0$. The presence of
$\lambda+1$ in the denominator should come as no surprise, minding
that $\lambda$ is the eigen-value of Hilbert-Schimdt operator. The
Fredholm alternative gives us a way of solving the integral
equation in terms of eigen-functions. Since $|\lambda|\leq\lambda_{1}=0.25553210...<1$, the
integral equation is {\it a posteriori} solvable, and $\lambda+1$
appears in the denominators. Curiously, it is possible to approach
this identity numerically. One of the motivations is to check its
validity, since the result heavily depends on the validity
of almost all the preceding results in \cite{ga2}. The left integral causes no
problems, since Taylor coefficients of $G_{\lambda}(z)$ can be
obtained at high precision as an eigen-vector of a finite matrix,
which is the truncation of an infinite one. On the other hand, the
right integral can be evaluated with less precision, since it
involves $F(x)$, and thus requires more time and space consuming
continued fractions algorithm. Nevertheless, the author of this
paper have checked it with completely satisfactory outcome,
confirming the validity.\\

Just as interestingly, results $(i)$ and $(iv)$ can be though as a reflection of a ``pair-correlation" between eigen-values $\lambda$ and eigen-value $-1$ (see Section 3 for some remarks on this topic). Moreover, minding properties of distributions $F_{\mu}(x)$ (here $\mu$ simply means another eigen-value), the following result can be obtained. Given the conditions enforced on $F_{\mu}$ by (\ref{psi}), identity $(\ref{summa})$ is replaced by (rigid for $f\in{\sf C}^{\omega}$)
\begin{eqnarray*}
\int\limits_{0}^{1}f(x)\d F_{\mu}(x)=-\frac{1}{\mu}\sum\limits_{n=1}^{\infty}
\int\limits_{0}^{1}f\Big{(}\frac{1}{x+n}\Big{)}2^{-n}\d F_{\mu}(x).
\end{eqnarray*}
Then our trick works smoothly again, and this yields an identity
\begin{eqnarray*}
\int\limits_{0}^{1}G_{\lambda}(-x)(\lambda+\mu x^{2})\d F_{\mu}(x)=0.
\end{eqnarray*}
The fact, consequently, is an interesting example of ``pair correlation" between eigen-values of Hilbert-Schmidt operator in (\ref{fred}). Using definition of distribution $F_{\mu}$, the last identity is equivalent to
\begin{eqnarray*}
\sum\limits_{L=1}^{\infty}(-1)^{L}\big{(}m^{(\mu)}_{L}m^{(\lambda)}_{L+1}\lambda-
m^{(\lambda)}_{L}m^{(\mu)}_{L+1}\mu\big{)}=0,
\end{eqnarray*}
and thus is a property of "orthogonality" of $G_{\lambda}(z)$. This expression is symmetric regarding $\mu$ and $\lambda$. As could be expected, it is void in case $\mu=\lambda$. As a matter of fact, the proof of the above identity is fallacious, since the definition of distributions $F_{\lambda}$ does not imply properties (\ref{psi}) (these simply have no meaning). Nevertheless, numerical calculations show that the last identity truly holds. We hope to return to this topic in the future.
\section{Fourier series}
Minkowski question mark function $F(x)$, originally defined for $x\geq 0$ by (\ref{min}),
can be extended naturally to $\mathbb{R}$ simply by functional equation
$F(x+1)=\frac{1}{2}+\frac{1}{2}F(x)$. Such an extension still is given by the
expression (\ref{min}), with the difference that $a_{0}$ can be negative integer.
Naturally, the second functional equation is not preserved for negative $x$.
Thus, we have
\begin{eqnarray*}
2^{x+1}(F(x+1)-1)=2^{x}(F(x)-1)\text{ for }x\in\mathbb{R}.
\end{eqnarray*}
So, $2^{x}(F(x)-1)$ is a periodic function, which we will denote by
$-\Psi(x)$. Figure 3 gives the graph of $\Psi(x)$ for $x\in[0,1]$.
\begin{figure}[h]
\begin{center}
\includegraphics[width=190pt,height=440pt,angle=-90]{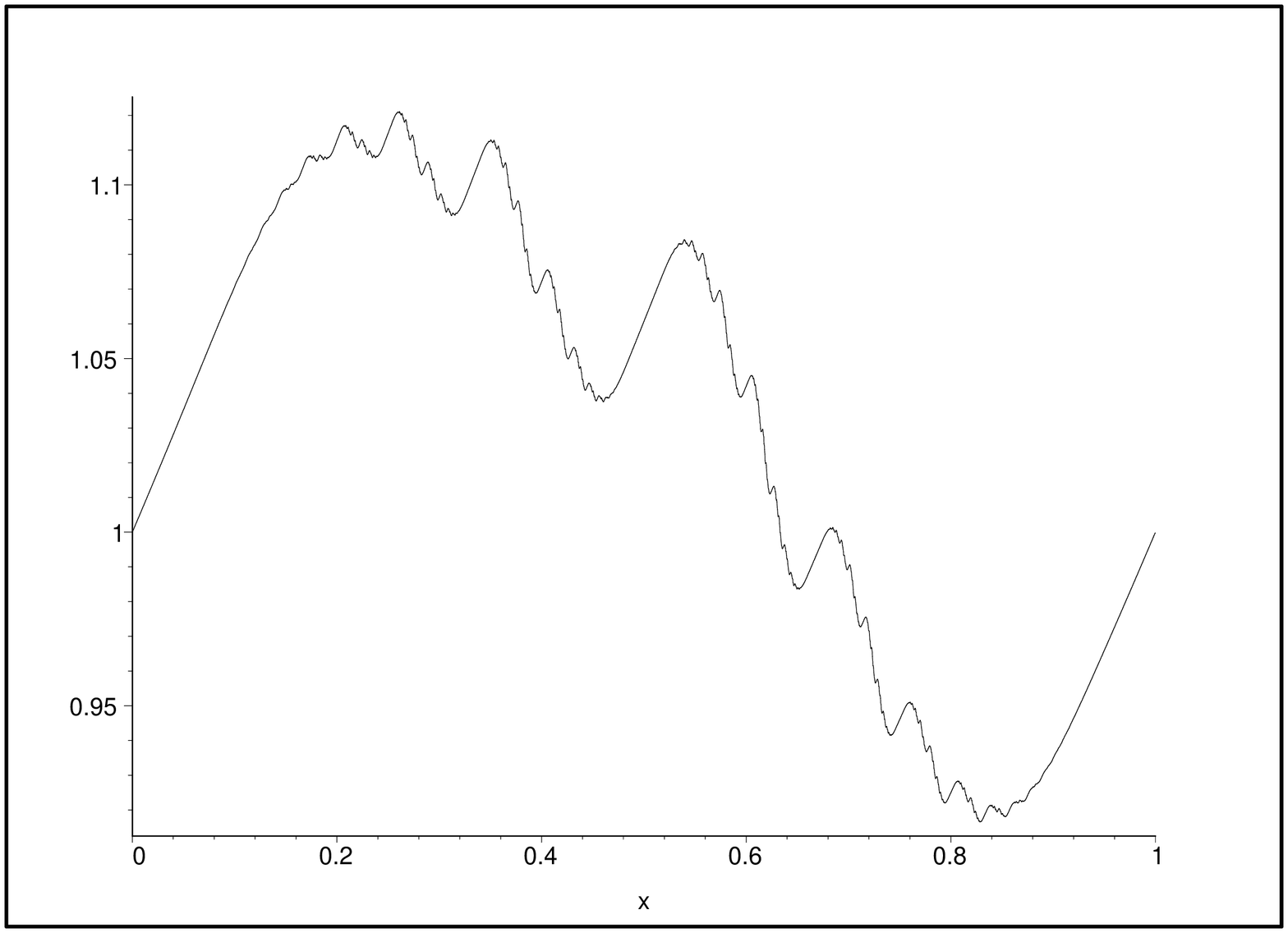}
\caption{$\Psi(x)$}
\end{center}
\end{figure}
Thus, $F(x)=-2^{-x}\Psi(x)+1$. Since $F(x)$ is singular, the same is true for
$\Psi(x)$: it is differentiable almost everywhere, and for these points
one has $\Psi'(x)=\log 2\cdot\Psi(x)$. As periodic function, it has an associated Fourier series expansion
$\Psi(x)\sim\sum\limits_{n=-\infty}^{\infty}c_{n}e^{2\pi i nx}.$ Since $F(x)$ is
real function, $c_{-n}=\overline{c_{n}}$, $n\in\mathbb{Z}$. Let for $n\geq 1$,
$c_{n}=a_{n}+ib_{n}$, and $a_{0}=\frac{c_{0}}{2}$. Here we list initial
numerical values for $c^{\star}_{n}=c_{n}(2\log 2-4\pi in)$ (see the next
proposition for the reason
of this normalization).\\

\begin{tabular}{l l l}
$c^{\star}_{0}=1.428159$,             & $c^{\star}_{3}=+0.128533-0.026840i$, & $c^{\star}_{6}=-0.262601+0.004128i$, \\
$c^{\star}_{1}=-0.521907+0.148754i$, & $c^{\star}_{4}=-0.140524-0.021886i$, & $c^{\star}_{7}=+0.198742-0.013703i$, \\
$c^{\star}_{2}=-0.334910-0.017869i$, & $c^{\star}_{5}=+0.285790+0.003744i$, &  $c^{\star}_{8}=-0.008479+0.024012i$. \\
\end{tabular}\\

It is important to note that we do not pose the question about the convergence of this Fourier series. For instance, the authors of \cite{salem} and \cite{reese} give examples of singular monotone increasing functions $f(x)$, whose Fourier-Stieltjes coefficients $\int_{0}^{1}e^{2\pi i nx}\d f(x)$ do not vanish, as $n\rightarrow\infty$. In \cite{salem}, the author even investigated $f(x)=?(x)$. In our case, the convergence problem is far from clear. Nevertheless, in all cases we substitute $-2^{-x}\Psi(x)$ instead of $(F(x)-1)$ under integral. Let, for example, $W(x)$ be a continuous function of at most polynomial growth, as $x\rightarrow\infty$, and let $\Psi_{N}(x)=\sum\limits_{n=-N}^{N}c_{n}e^{2\pi inx}$. Then
\begin{eqnarray*}
\Big{|}\int\limits_{0}^{\infty}W(x)\Big{(}(F(x)-1)+2^{-x}\Psi_{N}(x)\Big{)}\d x\Big{|}\ll
\sum\limits_{r=0}^{\infty}|W(r)|2^{-r}\cdot\int\limits_{0}^{1}|2^{x}(F(x)-1)+\Psi_{N}(x)|\d x.
\end{eqnarray*}
Since $2^{x}(F(x)-1)\in\mathcal{L}_{2}[0,1]$, the last integral tends to $0$, as $N\rightarrow\infty$. As it was said, this makes the change of $(F(x)-1)$ into $-2^{-x}\Psi(x)$ under integral legitimate, and this also justifies term-by-term integration. Henceforth, we will omit a step of changing $\Psi(x)$ into $\Psi_{N}(x)$, and taking a limit $N\rightarrow\infty$.\\
A general formula for the Fourier coefficients is given by
\begin{prop} Fourier coefficients $c_{n}$ are related to special values of
exponential generating function $m(t)$ through the equality
\begin{eqnarray*}
c_{n}=\frac{m(\log 2-2\pi in)}{2\log 2-4\pi in}, \text{ and }c_{n}=O(n^{-1}).
\end{eqnarray*}
\end{prop}
\proof We have (note that $F(1)=\frac{1}{2}$):
\begin{eqnarray*}
c_{n}=-\int\limits_{0}^{1}2^{x}(F(x)-1)e^{-2\pi inx}\d x= -\frac{1}{\log 2-2\pi
in}\int\limits_{0}^{1}(F(x)-1)\d e^{x(\log 2-2\pi in)}=\\
\frac{1}{\log 2-2\pi in}\int\limits_{0}^{1}e^{x(\log 2-2\pi in)}\d F(x)=
\frac{m(\log 2-2\pi in)}{2\log 2-4\pi in}.
\end{eqnarray*}
The last assertion of proposition is obvious.$\blacksquare$\\
This proposition is good example of intrinsic relations among the three
functions $F(x)$, $G(z)$ and $m(t)$. Indeed, the moments $m_{L}$ of $F(x)$ give Taylor
coefficients of $G(z)$, which are proportional (up to the factorial multiplier) to
Taylor coefficients of $m(t)$. Finally, special values of $m(t)$ on a discrete
set of vertical line produce ``Fourier coefficients" of $F(x)$.\\
Next proposition describes explicit relations among Fourier coefficients and the moments. Additionally, in the course of the proof we obtain the expansion of $G(z)$ for negative real $z$ in terms of incomplete gamma integrals.
\begin{prop} For $L\geq 1$, one has
\begin{eqnarray}
M_{L}=L!\sum\limits_{n\in\mathbb{Z}}\frac{c_{n}}{(\log 2-2\pi in)^{L}}.
\label{sum}
\end{eqnarray}
\end{prop}
\proof Let $z<0$ be fixed negative real. Then integration by parts gives
\begin{eqnarray*}
G(z+1)=\int\limits_{0}^{\infty}\frac{x}{1-xz}\d (F(x)-1)=
\int\limits_{0}^{\infty}\frac{1}{(1-xz)^{2}}2^{-x}\Psi(x)\d x=\\
\sum\limits_{n=-\infty}^{\infty}c_{n}\int\limits_{0}^{\infty}\frac{1}{(1-xz)^{2}}
2^{-x}e^{2\pi inx}\d x=\sum\limits_{n=-\infty}^{\infty}c_{n}V_{n}(z),
\end{eqnarray*}
where
\begin{eqnarray*}
V_{n}(z)=\int\limits_{0}^{\infty}\frac{1}{(1-xz)^{2}}e^{-x(\log 2-2\pi in)}\d
x= \frac{1}{\log 2-2\pi in}\int\limits_{0}^{\infty(\log 2-2\pi
in)}\frac{1}{(1-\frac{yz}{\log 2-2\pi in})^{2}}e^{-y}\d y.
\end{eqnarray*}
Since by our convention $z<0$, the function under integral does not have poles for $\Re y>0$, and Jordan's lemma gives
\begin{eqnarray*}
V_{n}(z)=\frac{1}{\log 2-2\pi
in}\int\limits_{0}^{\infty}\frac{1}{(1-\frac{yz}{\log
2-2\pi in})}e^{-y}\d y=\\
\frac{1}{\log 2-2\pi in}\cdot V\Big{(}\frac{z}{\log 2-2\pi in}\Big{)},\text{
where }V(z)=\int\limits_{0}^{\infty}\frac{1}{(1-yz)^{2}}e^{-y}\d y.
\end{eqnarray*}
The function $V(z)$ is defined for the same values of $z$ as $G(z+1)$ and therefore is defined
in the cut plane $\mathbb{C}\setminus(0,\infty)$. Consequently, this implies
\begin{eqnarray}
G(z+1)=\sum\limits_{n\in\mathbb{Z}}\frac{c_{n}}{\log 2-2\pi in}\cdot
V\Big{(}\frac{z}{\log 2-2\pi in}\Big{)}. \label{period}
\end{eqnarray}
The formula is only valid for real $z<0$. The obtained series converges uniformly, since $|1-y\frac{z}{\log 2-2\pi in}|\geq 1$ for $n\in\mathbb{Z}$ and $z<0$. Since
\begin{eqnarray*}
V\Big{(}\frac{1}{z}\Big{)}=-ze^{-z}\int\limits_{1}^{\infty}\frac{1}{y^{2}}e^{yz}\d y,
\end{eqnarray*}
this gives us the expansion of $G(z+1)$ on negative real line in terms of incomplete gamma integrals. As noted before, and this can be seen from (\ref{est}), the function $G(z)$ has all left derivatives at $z=1$. Further,
$(L-1)-$fold differentiation of $V(z)$ gives
\begin{eqnarray*}
V^{(L-1)}(z)=L!\int\limits_{0}^{\infty}\frac{y^{L-1}}{(1-yz)^{L+1}}e^{-y}\d
y\Rightarrow V^{(L-1)}(0)=L!(L-1)!.
\end{eqnarray*}
Comparing with (\ref{period}) and (\ref{est}), this gives the desired relation among moments
$M_{L}$ and Fourier coefficients, as stated in the proposition. $\blacksquare$

It is important
to compare this expression with the first equality of (\ref{symm}). Indeed,
since $m(t)$ is entire, that equality via Cauchy residue formula implies the
result obtained in \cite{as}, i.e.
\begin{eqnarray}
M_{L}\sim \frac{m(\log 2)}{2\log 2}\Big{(}\frac{1}{\log 2}\Big{)}^{L}L!\label{fur}
\end{eqnarray}
It is exactly the leading term in (\ref{sum}), corresponding to $n=0$. \ignore{More
than that, we have
\begin{prop}
Let $|z|<\log 2$ be complex number. Then
\begin{eqnarray*}
M(z)=\sum\limits_{L=0}^{\infty}\frac{z^{L}}{L!}M_{L}=
\frac{1}{2}\sum\limits_{n\in\mathbb{Z}}\frac{m(\log 2+2\pi in)}{\log 2+2\pi
in-z}.
\end{eqnarray*}
\end{prop}}
\section{Associated zeta function}
 Recall that for complex $c$ and $s$, $c^{s}$ is multi-valued
complex function, defined as $e^{s\log c}=e^{s(\log|c|+i\arg(c))}$. Henceforth, we fix the branch of the logarithm by requiring that the value of $\arg c$
for $c$ in the right half plane $\Re c>0$ is in the
range $(-\pi/2,\pi/2)$. Thus, if
$s=\sigma+it$, and if we denote $r_{n}=\log 2+2\pi i n$, then
$|r_{n}^{-s}|=|r_{n}|^{-\sigma}e^{t\arg r_{n}}\sim|r_{n}|^{-\sigma}e^{\pm \pi
t/2}$ as $n\rightarrow\pm\infty$. Minding this convention and the identity (\ref{sum}), we
introduce the zeta
function, associated with Minkowski question mark function.\\
\begin{defin} The dyadic zeta function $\zm(s)$ is defined in the half plane $\Re
s>0$ by the series
\begin{eqnarray}
\zm(s)=\sum\limits_{n\in\mathbb{Z}}\frac{c_{n}}{(\log 2-2\pi in)^{s}},
\label{dir}
\end{eqnarray}
where $c_{n}$ are Fourier coefficients of $\Psi(x)$, and for each $n$, $(\log 2-2\pi in)^{s}$ is understood in the meaning just
described.
\end{defin}
Then we have
\begin{thm}
$\zm(s)$ has an analytic continuation as an entire function to the whole plane
$\mathbb{C}$, and satisfies the functional equation
\begin{eqnarray}
\zm(s)\Gamma(s)=-\zm(-s)\Gamma(-s).
\label{tet}
\end{eqnarray}
Further, $\zm(L)=\frac{M_{L}}{L!}$ for $L\geq 1$. $\zm(s)$ has trivial zeros
for negative integers: $\zm(-L)=0$ for $L\geq 1$, and $\zm'(-L)=(L-1)!(-1)^{L}M_{L}$.
Additionally, $\zm(s)$ is real on the real line, and thus
$\zm(\overline{s})=\overline{\zm(s})$. The behavior of $\zm(s)$ in the vertical
strips is given by estimate
\begin{eqnarray*}
|\zm(\sigma+it)|\ll t^{-\sigma-1/2}\cdot e^{\pi|t|/2}
\end{eqnarray*}
uniformly for $a\leq \sigma\leq b$, $|t|\rightarrow\infty$.
\end{thm}
As we will see, these properties are immediate (subject to certain regularity conditions) for any distribution $f(x)$
with a symmetry property $f(x)+f(1/x)=1$. Nevertheless, it is
a unique characteristic of $F(x)$ that the corresponding zeta function can be
given a Dirichlet series expansion, like (\ref{dir}). We do not give the proof of the converse result, since there is no motivation for this. But empirically, we see that this functional equation is equivalent exactly to the symmetry property.
Additionally, the presence of a Dirichlet series expansion yields a functional equation of the kind
$f(x+1)=\frac{1}{2}f(x)+\frac{1}{2}$. Generally speaking, these two together are unique for $F(x)$. Note also that the functional equation implies that $\zm(it)\Gamma(1+it)=\int\limits_{0}^{\infty}x^{it}\d F(x)$ is real for real $t$. Figure 4 shows its graph for $1.5\leq t\leq 90$. Further calculations support the claim that this function has infinitely many zeros on the critical line $\Re s=0$. On the other hand, numerical calculations of the contour integral reveal that there exists much more zeros apart from these.\\
\begin{figure}[h]
\begin{center}
\includegraphics[width=190pt,height=440pt,angle=-90]{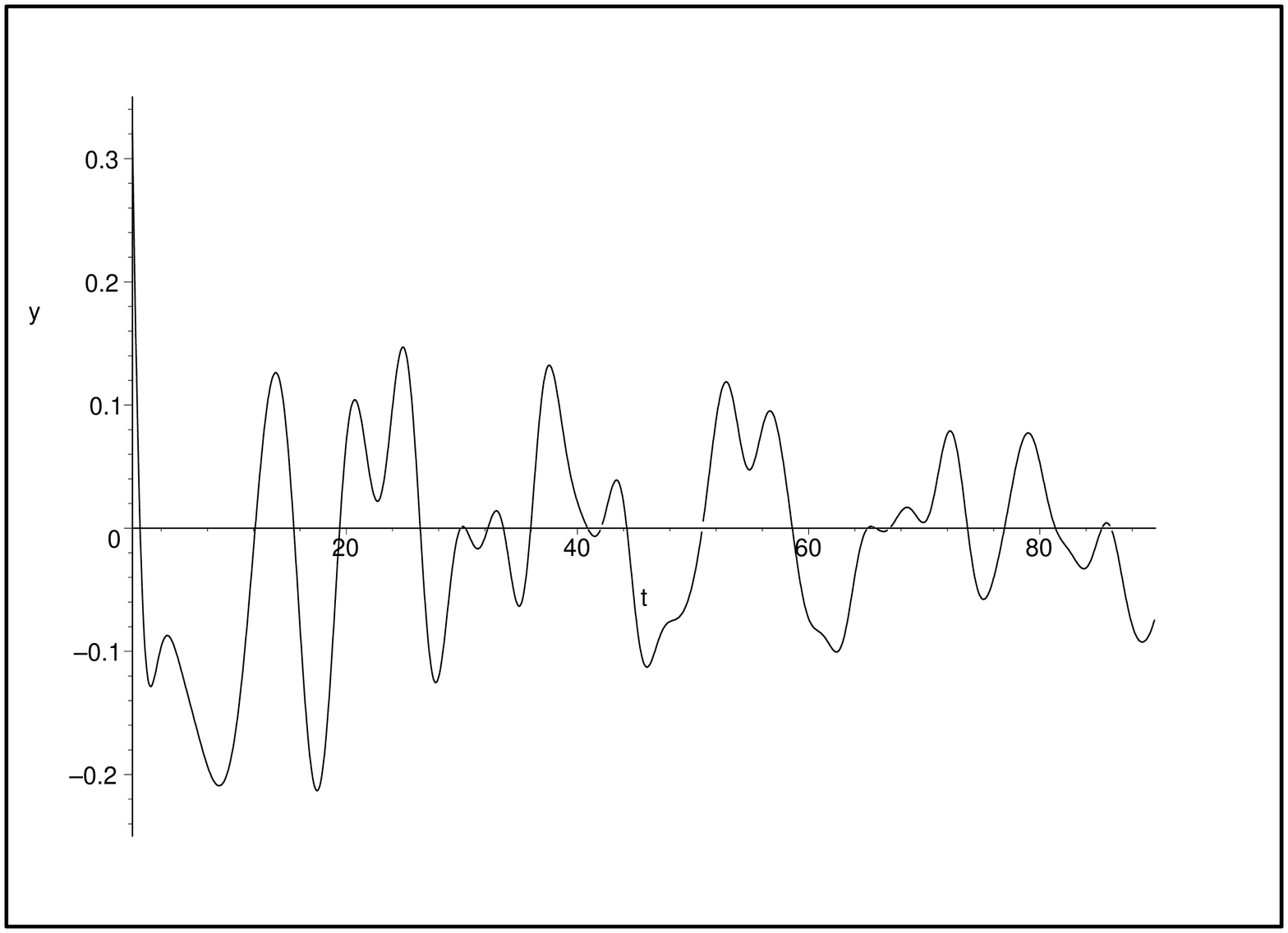}
\caption{$\zm(it)\Gamma(1+it)$}
\end{center}
\end{figure}
 \\
We need one classical integral.
\begin{lem}
Let $A$ be real number, $\arctan(A)=\phi\in(-\frac{\pi}{2},\frac{\pi}{2})$, and $\Re s>0$. Then
\begin{eqnarray*}
\int\limits_{0}^{\infty}x^{s-1}e^{-x}\cos(Ax)\d x=
\frac{1}{(1+A^{2})^{s/2}}\cos(\phi s)\Gamma(s).
\end{eqnarray*}
The same is valid with $\cos$ replaced by $\sin$ on both sides.
\end{lem}
This can be found in any extensive table of gamma integrals or tables of Mellin transforms.\\

{\it Proof of Theorem. }Let for $n\geq0$, $\arctan\frac{2\pi n}{\log 2}=\phi_{n}$. We will
calculate the following integral. Let $\Re s>0$. Then integrating by parts and
using lemma, one obtains
\begin{eqnarray*}
\int\limits_{0}^{\infty}x^{s}\d (F(x)-1)=s\int\limits_{0}^{\infty}x^{s-1}2^{-x}\Psi(x)\d x=
s\sum\limits_{n\in\mathbb{Z}}c_{n}\int\limits_{0}^{\infty}x^{s-1}2^{-x}e^{2\pi i nx}\d x\\
=s\sum\limits_{n=0}^{\infty}\int\limits_{0}^{\infty}x^{s-1}\Big{(}2a_{n}\cos(2\pi nx)-
2b_{n}\sin(2\pi nx)\Big{)}2^{-x}\d x=\\
2s\Gamma(s)\sum\limits_{n=0}^{\infty}|\log 2+2\pi
ni|^{-s}\Big{(}a_{n}\cos(\phi_{n} s)- b_{n}\sin(\phi_{n}
s)\Big{)}=s\Gamma(s)\sum\limits_{n\in\mathbb{Z}}\frac{c_{n}}{(\log 2-2\pi i
n)^{s}}.
\end{eqnarray*}
Note that the function
$\int_{0}^{\infty}x^{s}\d F(x)$ is clearly analytic and entire. Thus, $s\Gamma(s)\zm(s)$
is an entire function, and this proves the first statement of the theorem. Since $F(x)+F(1/x)=1$, this gives
$\int_{0}^{\infty}x^{s}\d F(x)=\int_{0}^{\infty}x^{-s}\d F(x)$,
and this, in turn, implies functional equation. All other statements follow easily from this,
our previous results, and known properties of the $\Gamma-$function. In particular, if $s=\sigma+it$,
\begin{eqnarray*}
|\zm(s)\Gamma(s+1)|\leq\int\limits_{0}^{\infty}|x^{s}|\d F(x)=\zm(\sigma)\Gamma(\sigma+1),
\end{eqnarray*}
and the last statement of the theorem follows from the Stirling's formula for $\Gamma-$function:
$|\Gamma(\sigma+it)|\sim\sqrt{2\pi}t^{\sigma-1/2}e^{-\pi|t|/2}$ uniformly for $a\leq \sigma\leq b$,
as $|t|\rightarrow\infty$. $\blacksquare$\\

At this stage, we will make some remarks, concerning the analogy and
differences with the classical results known for the Riemann zeta function
$\zeta(s)=\sum\limits_{n=1}^{\infty}\frac{1}{n^{s}}$. Let $\theta(x)$ denote
the usual theta function: $ \theta(x)=\sum\limits_{n\in\mathbb{N}}e^{\pi i
n^{2}x}$, $\Im x>0$. The following table summarizes all ingredients,
which eventually produce the functional equation both for $\zeta(s)$ and
$\zm(s)$.\\

\it
\begin{tabular}{|r | r | r|}
\hline
Function & $\zeta(s)$& $\zm(s)$\\
\hline
Dirichlet series expansion& Periodicity: $\theta(x+2)=\theta$ & Periodicity: $F'(x+1)=\frac{1}{2}F'(x)$\\
Functional equation & $\theta(ix)=\frac{1}{\sqrt{x}}\theta(\frac{i}{x})$& $F'(x)=-F'(\frac{1}{x})$\\
\hline
\end{tabular}\\
\rm

Since $F(x)$ is a singular function, its derivative should be considered as a
distribution on the real line. For this purpose, it is sufficient to consider
a distribution $U(x)$ as a derivative of a continuous function $V(x)$, for which the
scalar product $\langle U,f\rangle$, defined for functions $f\in
C^{\infty}(\mathbb{R})$ with compact support, equals to $-\langle
V,f'\rangle=-\int\limits_{\mathbb{R}}f'(x)V(x)\d x$. Thus, both $\theta(x)$ and
$2^{x}F'(x)$ are periodic distributions. This guarantees that the appropriate
Mellin transform can be factored into the product of Dirichlet series and gamma
factors. Finally, the functional equation for the distribution produces the
functional equation for the Mellin transform. The difference arises from the
fact that for $\theta(x)$, the functional equation is symmetry property on the
imaginary line, whereas for $F'(x)$ we have the symmetry on the real line instead.
This explains the unusual fact that in (\ref{dir}) we have the
summation over the discrete set of the vertical line, instead of the
summation over integers.\\

We will finish by proving another result, which links $\zm(s)$ to the Mellin
transform of $G(-z+1)$. This can be done using expansion (\ref{period}), but we rather chose a direct way. Let
$\int\limits_{0}^{\infty}G(-z+1)z^{s-1}\d z=G^{*}(s)$. Symmetry property of Theorem 1 implies
that $G(-z+1)$ has a simple zero, as $z\rightarrow\infty$ along the positive real line. Thus, basic properties of Mellin transform imply that $G^{*}(s)$ is defined for $0< \Re s< 1$. For these values of $s$, we have the
following classical integral:
\begin{eqnarray*}
\int_{0}^{\infty}\frac{z^{s-1}}{1+z}\d z\mathop{=}^{\frac{z}{1+z}\rightarrow x}\int\limits_{0}^{1}x^{s-1}(1-x)^{-s}\d x=\Gamma(s)\Gamma(1-s)=\frac{\pi}{\sin\pi s}.
\end{eqnarray*}
 Thus, using (\ref{gen}), we get
\begin{eqnarray*}
G^{*}(s)=\int\limits_{0}^{\infty}\int\limits_{0}^{\infty}\frac{xz^{s-1}}{1+xz}\d F(x)\d z=\int\limits_{0}^{\infty}\int\limits_{0}^{\infty}\frac{z^{s-1}}{1+z}x^{1-s}\d z\d F(x)=\frac{\pi}{\sin \pi s}\int\limits_{0}^{\infty}x^{1-s}\d F(x).
\end{eqnarray*}
This holds for $0< \Re s< 1$. Due to the analytic continuation, this gives
\begin{prop}
For $s\in\mathbb{C}\setminus\mathbb{Z}$, we have an identity
\begin{eqnarray*}
G^{*}(s)=\zm(s-1)\Gamma(s)\cdot\frac{\pi}{\sin\pi s}.
\end{eqnarray*}
\end{prop}
Therefore, $G^{*}(s)$ is a meromorphic function, $G^{*}(s+1)=-G^{*}(-s+1)$, and $\text{res}_{s=L}G^{*}(s)=(-1)^{L}M_{L-1}$. This is, of course, the general property of the Mellin transform, since formally $G(z+1)=\sum_{L=0}^{\infty}M_{L}z^{L-1}$. Thus, $G(z+1)\sim\sum\limits_{L=0}^{M}M_{L}z^{L-1}$ in the left neighborhood of $z=0$.

\section{Concluding remarks}
\subsection{Dyadic period functions in $\mathbb{H}$}
As noted in \cite{ga2}, one encounters the surprising fact that in the
upper half plane $\mathbb{H}$, the equation (\ref{funct}) is also satisfied by
$\frac{i}{2\pi}G_{1}(z)$, where $G_{1}(z)$ stands for the Eisenstein series, which in $\mathbb{H}$ is given by Fourier expansion $\frac{\pi^{2}}{3}-8\pi^{2}\sum_{n=1}^{\infty}\sigma_{1}(n)e^{2\pi inz}$. The reason for this is that $G_{1}(z)$ transforms under the action of ${\sf PSL}_{2}(\mathbb{Z})$ as (see \cite{serre})
\begin{eqnarray*}
G_{1}(z+1)=G_{1}(z),\quad G_{1}(-1/z)=z^{2}G_{1}(z)-2\pi iz.
\end{eqnarray*}
Let $f_{0}(z)=G(z)-\frac{i}{2\pi}G_{1}(z)$, where $G(z)$ is the function in Theorem 1. Then
for $z\in\mathbb{H}$, $f_{0}(z)$ satisfies the homogeneous form of the three term functional equation
(\ref{funct}); moreover $f_{0}(z)$ is bounded, when $\Im z\rightarrow\infty$. Thus, if $f(z)=f_{0}(z)$,
\begin{eqnarray*}
-\frac{1}{(1-z)^{2}}f\Big{(}\frac{1}{1-z}\Big{)}+2f(z+1)=f(z).
\label{hom}
\end{eqnarray*}
 Therefore, denote by ${\tt DPF}^{0}$ the $\mathbb{C}-$linear vector space of solutions of this three term functional equation, which are
 holomorphic in $\mathbb{H}$ and are bounded at infinity, and call it {\it the space of dyadic period functions in the upper half-plane}. Consequently, this space is at least one-dimensional. If we abandon the growth condition, then the corresponding space ${\tt DPF}$ is infinite-dimensional.  This is already true for periodic solutions. Indeed, if $f(z)$ is a periodic solution, then $f(z)=\frac{1}{z^{2}}f(-1/z)$. Let $P(z)\in\mathbb{C}[z]$, and suppose that $j(z)$ stands, as usually, for the $j-$invariant. Then any modular function of the form $j'(z)P(j(z))$ satisfies this equation. Additionally, there are non-periodic solutions, given by $f_{0}(z)P(j(z))$. Therefore, $G(z)$ surprisingly enters the profound domain of classical modular forms and functions for ${\sf PSL}_{2}(\mathbb{Z})$. Moreover, in the space ${\tt DPF}$, one establishes the relation between real quadratic irrationals (via $G(z)$, Minkowski question mark function $F(x)$ and continued fraction algorithm), and imaginary quadratic irrationals (via $j-$invariant and its special values). Hence, it is greatly desirable to give the full description and structure of spaces ${\tt DPF}^{0}$ and ${\tt DPF}$.
\subsection{Where should the true arithmetic zeta function come from?}
Here we present some remarks, concerning the zeta function $\zm(s)$. This object is natural for the question mark function - its Dirichlet coefficients are the Fourier coefficients of $F(x)$, and its special values at integers are proportional to the moments $M_{L}$. Moreover, its relation to $G(z)$, $m(t)$ and $F(x)$ is the same as the role of $L-$series of Maass wave forms against analogous objects \cite{zagier}. Nevertheless, one expects richer arithmetic object associated with Calkin-Wilf tree, since the latter consists or rational numbers, and therefore can be canonically embedded into the group of \'{i}deles $\mathbb{A}_{\mathbb{Q}}$. The $p-$adic distribution of rationals in the $n-$th generation of Calkin-Wilf tree was investigated in \cite{as}. Surprisingly, Eisenstein series $G_{1}(z)$ yet again manifest, as in case of $\mathbb{R}$ (see previous subsection). Nevertheless, there is no direct way of normalizing moments of the $n-$th generation in order for them to converge in the $p-$adic norm. There is an exception. As shown in \cite{as},
\begin{eqnarray*}
\sum\limits_{a_{0}+a_{1}+...+a_{s}=n}[a_{0},a_{1},..,a_{s}]=
3\cdot2^{n-2}-\frac{1}{2},
\end{eqnarray*}
and thus we have a convergence only in the $2-$adic topology, namely to the value $-\frac{1}{2}$. The investigation of $p-$adic values of moments is relevant for the following reason. Let us apply $F(x)$ to each rational number in the Calkin-Wilf tree. What we obtain is the following:
$$
\xymatrix @R=.5pc @C=.5pc { & & & & & & & {1\over 2} & & & & & & & \\
& & & {1\over 4} \ar@{-}[urrrr] & & & & & & & & {3\over 4} \ar@{-}[ullll] & & & \\
& {1\over 8} \ar@{-}[urr] & & & & {5\over 8}\ar@{-}[ull] & & & & {3\over 8}\ar@{-}[urr] & & & & {7\over 8} \ar@{-}[ull] & \\
{1\over 16} \ar@{-}[ur] & & {9\over 16} \ar@{-}[ul] & & {5\over
16} \ar@{-}[ur] & & {13\over 16} \ar@{-}[ul] & & {3\over 16}
\ar@{-}[ur] & & {11\over 16} \ar@{-}[ul] & & {7\over 16}
\ar@{-}[ur] & & {15\over 16} \ar@{-}[ul] }
$$
Using (\ref{distr}), we deduce that this tree starts from the root $\frac{1}{2}$, and then inductively each rational $r$ produces two offsprings: $\frac{r}{2}$ and $\frac{r}{2}+\frac{1}{2}$. One is therefore led to the following\\

\textbf{Task. }\it Produce a natural algorithm, which takes into account $p-$adic and real properties of the above tree, and generates Riemann zeta function $\zeta(s)$.\rm\\

We emphasize that the choice of $\zeta(s)$ is not accidental. In fact, $\mathbb{R}-$distribution of the above tree is a uniform one with support $[0,1]$. Further, there is a natural algorithm to produce ``characteristic function of ring of integers of $\mathbb{R}"$ (that is, $e^{-\pi x^{2}}$) from the uniform distribution via the central limit theorem through the expression
\begin{eqnarray*}
\int\limits_{\mathbb{R}}f(x)e^{-\pi x^{2}}\d x=\lim\limits_{N\rightarrow\infty}\frac{1}{2^{N}}\int\limits_{-1}^{1}\d x_{1}...\int\limits_{-1}^{1}\d x_{N}f\Big{(}\frac{x_{1}+...+x_{N}}{\sqrt{\frac{2}{3}\pi N}}\Big{)}.
\end{eqnarray*}
(For clarity, here we take the uniform distribution in the interval $[-1,1]$). This formula and this explanation and treatment of $e^{-\pi x^{2}}$ as ``characteristic function of the ring of integer of $\mathbb{R}$" is borrowed from \cite{harad}, p. 7. Further, the operator which is invariant under uniform measure has the form $[\mathcal{U}f](x)=\frac{1}{2}f\big{(}\frac{x}{2}\big{)}+
\frac{1}{2}f\big{(}\frac{x}{2}+\frac{1}{2}\big{)}$. Indeed, for every $f\in{\sf C}[0,1]$, one has $\int_{0}^{1}[\mathcal{U}f](x)\d x=\int_{0}^{1}f(x)\d x$. The spectral analysis of $\mathcal{U}$ shows that its eigenvalues are $2^{-n}$, $n\geq 0$, with corresponding eigen-functions being Bernoulli polynomials $B_{n}(x)$ \cite{flaj}. These, as is well known from the time of Euler, are intricately related with $\zeta(s)$. Moreover, the partial moments of the above tree can be defined as $\sum_{i=1}^{2^{N}}\big{(}\frac{2i-1}{2^{N}}\big{)}^{L}$. These values are also expressed in term of Bernoulli polynomials. As we know, there are famous Kummer congruences among Bernoulli numbers, which later led to the introduction of the $p-$adic zeta function $\zeta_{p}(s)$. Thus, the real distribution of the above tree and its spectral decomposition is deeply related to the $p$-adic properties. This justifies the choice in the task of  $\zeta(s)$.\\
Therefore, returning to Calkin-Wilf tree, one expects that moments can be $p-$adically interpolated, and some natural arithmetic zeta function can be introduced, as a ``pre-image" of $\zeta(s)$ under map $F$.\\
\ignore{\subsection{Dyadic invariant}
Now, for $z\in\mathbb{H}\cup\mathbb{R}$, define {\it holomorphic question mark
function } $F_{+}(z)$ by the expression
\begin{eqnarray*}
F_{+}(z)=-2^{-z}\Big{(}\frac{c_{0}}{2}+\sum\limits_{n=1}^{\infty}c_{n}e^{2\pi
inz}\Big{)}+\frac{1}{2}.
\end{eqnarray*}
Then $F_{+}(\overline{z})$ is defined in $\mathbb{H}_{-}\cup\mathbb{R}$, and in
their common set of definition (in $\mathbb{R}$) one has
$F_{+}(z)+\overline{F_{+}(\overline{z})}=F(z)$. Since $F_{+}(z)$ is not
homogeneous regarding the translation $z\rightarrow z+1$, the following object
will be one of the main topics in the consequent paper by the author. Let us
define {\it Dyadic invariant } $\mathcal{F}(z)$ for $z\in\mathbb{H}$ by the
expression $\mathcal{F}(z)=\frac{\d}{\d z}F_{+}(z)$. Then, obviously,
$\mathcal{F}(z+1)=\frac{1}{2}\mathcal{F}(z)$, $z\in\mathbb{H}$. We indulged
into liberty calling this function dyadic invariant, since it should inherit
some information regarding the transformation $z\rightarrow\frac{z}{z+1}$,
which is an intrinsic property of $F(x)$ and is given by (\ref{distr}).}

\par\bigskip

\noindent Giedrius Alkauskas, PhD student of School of Mathematical Sciences,
University of Nottingham, University Park, Nottingham NG7 2RD
United Kingdom\\
 {\tt giedrius.alkauskas@maths.nottingham.ac.uk}
\smallskip
\end{document}